\documentclass[11pt]{article}

\scrollmode
\usepackage{amsmath,amssymb,amsfonts,graphicx}
\usepackage[latin1]{inputenc}

\newtheorem{thm}{Theorem}[section]
\newtheorem{prop}[thm]{Proposition}
\newtheorem{lem}[thm]{Lemma}
\newtheorem{df}[thm]{Definition}

\newtheorem{rem}[thm]{Remark}

\newtheorem{ass}[thm]{Assumption}
\newtheorem{cor}[thm]{Corollary}

\def\be#1 {\begin{equation} \label{#1}}
\newcommand{\ee}{\end{equation}}
\def\dem {\noindent {\bf Proof : }}

\newcommand{\mb}{\medskip\noindent}
\newcommand{\gb}{\bigskip\noindent}
\newcommand{\R}{\mathbb R}

\newcommand{\F}{\mathcal F}

\newcommand{\p}{\mu}

\newcommand{\D}{\mathcal D}

\newcommand{\diam}{diam}
\newcommand{\aver}[1]{-\hskip-0.46cm\int_{#1}}
\newcommand{\avert}[1]{-\hskip-0.38cm\int_{#1}}
\newcommand{\ind}{ {\bf 1} }
\DeclareMathOperator{\supp}{supp}
\DeclareMathOperator{\Lip}{Lip}

\newcommand{\B}{\beta}
\def\sqw{\hbox{\rlap{\leavevmode\raise.3ex\hbox{$\sqcap$}}$%
\sqcup$}}
\def\findem{\ifmmode\sqw\else{\ifhmode\unskip\fi\nobreak\hfil
\penalty50\hskip1em\null\nobreak\hfil\sqw
\parfillskip=0pt\finalhyphendemerits=0\endgraf}\fi}

\date{February 03, 2010}

\title{Abstract Hardy-Sobolev spaces and interpolation }
\date{February 03, 2010}

\begin {document}

\author{ N. Badr\\Institut Camille Jordan\\Universit\'e Claude Bernard Lyon 1 \\ 43 boulevard du 11 Novembre 1918\\ F-69622 Villeurbanne Cedex\\ badr@math.univ-lyon1.fr \and
F. Bernicot\\ Universit\'e de Paris-Sud \\F-91405 Orsay
Cedex\\frederic.bernicot@math.u-psud.fr \\}
\maketitle

\begin{abstract}
The purpose of this work is to describe an abstract theory of Hardy-Sobolev spaces on doubling Riemannian manifolds via
an  {\em atomic decomposition}. We study the real interpolation of these spaces with Sobolev spaces and finally give
applications to Riesz inequalities.
\end{abstract}

\mb {\bf Key-words :}  Sobolev spaces ; Hardy-Sobolev spaces ; Real interpolation ; Riesz inequalities. 

\tableofcontents

\gb The aim of the present work is to define  atomic Hardy-Sobolev spaces and interpolate them with Sobolev spaces on
Riemannian manifolds.
\\
One of the motivations is our Sobolev interpolation result \cite{Nadine}, \cite{Nadine1}  in different geometric frames, under the  doubling property and Poincar\'e inequalities. After this result, it is interesting to consider  a  ``nice" subspace of $W^{1,1}$ -- as is the Hardy space for $L^1$-- and study the interpolation of this  space with Sobolev spaces. Apart from the interpolation itself, the  use of the Hardy-Sobolev spaces that we construct gives strong boundedness of some linear operators  instead of the weak boundedness on $W^{1,1}$. For instance this is the case of the square root of the Laplace-Beltrami operator $ \Delta^{1/2}$.  \\
Another motivation responds to the recent improvements on the theory of Hardy spaces. In the last years, many works
were related to the study of specific Hardy spaces defined according to a particular operator (Riesz transforms,
Maximal regularity operator, Calder\'on-Zygmund operators, ... \cite{BJ,BJ2,DY,D1,FS,HM,R}). Mainly one of the most
interesting questions in this theory is the interpolation of these spaces with Lebesgue spaces in order to prove boundedness of some operators. \\
Although the theory of Hardy spaces is now well developed, the more recent theory of Hardy-Sobolev spaces is still not
unified.

\gb Before we state our results, let us briefly review the existing literature related to this subject.

\gb The Hardy-Sobolev spaces were studied by many authors in the Euclidean case. We mention R.~Strichartz \cite{strichartz2}. Related works
are \cite{ART}, \cite{KS}, \cite{cho}, \cite{miyachi}. They deal with ``classical'' Hardy-Sobolev spaces $ HS^1$ on
$\R^n$, which correspond to the Sobolev version of the Coifman-Weiss Hardy space $H^1_{CW}(\R^n)$~: $HS^{1}$ is the set
of functions $f \in H_{CW}^1$ such that each partial derivative of $f$ belongs to $H^1_{CW}$. Some of them consider the
homogeneous version of $HS^1$ and others only assume $f\in L^1$ instead of $f\in H^1_{CW}$.

\gb We recall that R.~Coifman proved an atomic decomposition for the classical Hardy space $H^1_{CW}$, which can be
defined by maximal functions (see \cite{FS}). In the Euclidean case, the question of atomic decomposition for the
homogeneous space $\dot{HS}^1$ was treated  in \cite{strichartz2} and \cite{cho}. However, in the non-Euclidean
case this issue is still not clear. In contrast, our idea is to introduce atomic Hardy-Sobolev spaces for which we can
prove real interpolation with Sobolev spaces. Then we are able to derive the interpolation of $HS^1$ with Sobolev
spaces.

\gb Let us now summarize the content of this paper. We refer the reader to the corresponding sections for definitions
and properties of the spaces and operators that we use in the statements.

\gb In the second part of Section 2, we define atomic Hardy-Sobolev spaces  $HS^1_{(\B),ato}$ for $1<\B\leq \infty$.
They correspond to the Sobolev version of the atomic Coifman-Weiss Hardy space $H^1_{CW}$ (defined by atomic
decomposition with $W^{1,\B}$-atoms). We compare these spaces  for different $\B$ in the following theorem:
\begin{thm} \label{thm:comp} Let $M$ be a complete Riemannian manifold satisfying $(D)$ and admitting a Poincar\'e inequality $(P_{q})$ for some $q>1$. Then
$HS_{(\B),ato}^{1}\subset HS_{(\infty),ato}^{1}$ for every $\B\geq q$ and therefore
$HS_{(\B_1),ato}^{1}=HS_{(\B_2),ato}^{1}$ for every $\B_1,\B_2\in[q,\infty]$.
 \end{thm}

\mb For the real interpolation of these spaces with Sobolev spaces, we obtain

\begin{thm} \label{thm:hardy} Let $M$ be a complete Riemannian manifold satisfying $(D)$ and $(P_q)$, for some $q \in(1,\infty)$. Let $r\in(1,\infty]$, $s\in(q,\infty]$,
$p\in(q,s)$ and $\theta\in(0,1)$ satisfying $\frac{1}{p}=(1-\theta)+\frac{\theta}{s}$. Then
$$
 W^{1,p}=\left(HS_{(r),ato}^{1}, W^{1,s}\right)_{\theta,p}= \left(HS^{1}, W^{1,s}\right)_{\theta,p}
 $$
 with equivalent norms.
 \end{thm}

\mb We also prove the homogeneous version of theses two theorems:
\begin{thm} \label{thm:comph} Let $M$ be a complete Riemannian manifold satisfying  $(D)$ and a Poincar\'e inequality $(P_{q})$ for some $q>1$. Then $  \dot{HS}_{(\B),ato}^{1}\subset  \dot{HS}_{(\infty),ato}^{1}$ for every $\B\geq q$ and therefore $ \dot{HS}_{(\B_1),ato}^{1}= \dot{HS}_{(\B_2),ato}^{1}$  for every $\B_1,\B_2\in[q,\infty]$.
 \end{thm}

\begin{thm} \label{thm:hardyhomo} Let $M$ be a complete Riemannian manifold satisfying $(D)$ and $(P_q)$, for some $1< q <\infty$. Let $r\in(1,\infty]$, $s\in(q,\infty]$ and $p\in(q,s)$ and $\theta\in(0,1)$  satisfying $\frac{1}{p}=(1-\theta)+\frac{\theta}{s}$. Then
$$
 \dot{W}^{1,p}=\left(\dot{HS}_{(r),ato}^{1}, \dot{W}^{1,s}\right)_{\theta,p}= \left(\dot{HS}^{1}, \dot{W}^{1,s}\right)_{\theta,p}
 $$
 with equivalent norms.
\end{thm}

\mb In the first part of section 2, given a collection of uniformly bounded operators on $ W^{1,\B}$:
$\mathbb{B}:=(B_Q)_{Q\in \mathcal{Q}}$ , we define  abstract atomic Hardy-Sobolev spaces $HW^{1}_{ato}$. For theses
spaces, we obtain in section 3 the following two interpolation results.

\begin{thm} \label{thm:generalbis} Let $M$ be a Riemannian manifold satisfying $(D)$. Let  $\sigma\in(1,\infty]$ and $p_0$  such that $\sigma'<p_0\leq \B$.
Let $\mathbb{B}:=(B_Q)_{Q\in \mathcal{Q}}$  be a collection of uniformly bounded operators on $ W^{1,\B}$ satisfying
\begin{equation}\label{assumpintro}
\frac{1}{\mu(Q)^{1/\sigma}}  \left\|f-B_Q^*(f)\right\|_{W^{-1,\sigma}(Q)} \lesssim M_{S,*,\B'}(f).
\end{equation}
Let $T$ be a bounded linear operator from $W^{1,p_0}$ to $L^{p_0}$ and from  $HW^{1}_{F,ato}$ to $L^1$. Then for every $p\in(\sigma',p_0)$ such that $(\B',p')\in {\mathcal I}_M$, there is a constant $c=c(p)$ such that for all function $f\in W^{1,p}\cap W^{1,p_0}$
$$ \|T(f)\|_{L^p} \leq c\|f\|_{W^{1,p}}.$$
Consequently, $T$ admits a continuous extension from $W^{1,p}$ to $L^p$.
\end{thm}
\begin{thm} \label{thm:general2} Let $M$ be a Riemannian manifold satisfying $(D)$ and of infinite measure $\mu(M)=\infty$.  Assume that the finite Hardy-Sobolev space is contained in $W^{1,1}$:
$$ HW^{1}_{F,ato} \hookrightarrow W^{1,1}$$ and that $\mathbb{B}$ satisfies (\ref{assumpintro}).
Let $\sigma\in(1,\infty]$ and $p_0$  satisfying $\sigma'<p_0\leq \B$.
Then for every $\theta\in(0,1)$ such that
$$ \frac{1}{p_\theta}:= (1-\theta) + \frac{\theta}{p_0}<\frac{1}{\sigma'}$$
and $(\B',p_\theta')\in {\mathcal I}_M$, we have
$$ \left(HW^{1}_{F,ato},W^{1,p_0}\right)_{\theta,p_\theta} = W^{1,p_\theta},$$
with equivalent norms.
\end{thm}

\mb Finally, the following theorem is an application of our result. It is proved in section 4 and applies to $
\Delta^{1/2}$.

 \begin{thm}\label{HSL}Let $M$ be a complete Riemannian manifold satisfying $(D)$.
 \begin{itemize}
 \item[1.] Assume that a Poincar\'e inequality $(P_1)$ holds. Let $T$ be a bounded  linear operator from $\dot{W}^{1,2}$ to $L^2$ and associated to a kernel satisfying
\begin{equation}\label{ker1}
 \sup_{Q \textrm{ ball}}\ \ \sup_{y,z\in Q} r_Q \int_{M\setminus 4Q} \left|K(x,y)-K(x,z) \right| d\mu(x)<\infty.
 \end{equation}
Then $T$ admits a unique extension from $\dot{HS}_{(2),ato}^1$ to $L^1$.
\item[2.] Assume that a Poincar\'e inequality $(P_2)$ holds. Let $T$ be a bounded linear operator from $W^{1,2}$ to $L^2$ and associated to a kernel satisfying (\ref{ker1}). \\
Then $T$ admits a unique extension from $HS_{(2),ato}^1$ to $L^1$.
\end{itemize}
\end{thm}
\begin{rem}
 Thanks to Theorem \ref{thm:comp}, in item 1. of Theorem \ref{HSL}, $T$ is then bounded from $\dot{HS}_{(\beta),ato}^1$ to $L^1$ for all $\beta\in (1,\infty]$. In item 2., $T$ is then bounded from $HS_{(\beta),ato}^1$ to $L^1$ for all $\beta\in [2,\infty]$.
\end{rem}

Consequently
\begin{cor} \label{cor1}
1-Let $T$ be as in item 1. of Theorem \ref{HSL}. Assume that a Poincar\'e inequality $(P_1)$ holds.
Then for all $p\in(1,2]$, the operator $T$ admits a continuous extension from $\dot{W}^{1,p}$ to $L^p$.\\
2-Let $T$ be as in item 2. of Theorem \ref{HSL}. Assume that a Poincar\'e inequality $(P_q)$ holds for some $q\in(1,2)$.
Then for all $p\in (q,2]$, the operator $T$ admits a continuous extension from $W^{1,p}$ to $L^p$.
\end{cor}

\mb We apply these last two theorems to the square root of the positive Laplace-Beltrami operator $\Delta^{1/2}$. In
\cite{AC}, P.~Auscher and  T.~Coulhon proved that under the  doubling property $(D)$ and a Poincar\'e inequality
$(P_q)$ for some $q\in[1,2)$, $(RR_p)$ (which is equivalent to the boundedness of $\Delta^{1/2}$ from $\dot{W}^{1,p}$
to $L^p$) holds for every $q<p\leq 2$. Moreover, $\Delta^{1/2}$ satisfies a weak type inequality $(RR_{qw})$ ($(RR_p)$ also holds in
this case for $2<p<\infty$).  Applying Theorem \ref{HSL}, we show that under $(D)$ and $(P_1)$ (resp. $(P_2)$ ) we have
a strong $(RR_1)$ (resp. $(nhRR_1)$) inequality for functions in the homogeneous (resp. non-homogeneous) atomic
Hardy-Sobolev space $\dot{HS}_{(\B),ato}^1$ (resp. $HS_{(\B),ato}^1$).

\gb We finish this introduction with a plan of the paper.  In section 1, we recall some definitions and properties that
we need. We define abstract Hardy-Sobolev spaces via atomic decomposition in the first part of section 2.  In the
second part we study  particular atomic Hardy-Sobolev spaces $HS^1_{(\B),ato}$ in more detail and  prove Theorem
\ref{thm:comp} . We also prove that under Poincar\'e inequality, these spaces are a particular case of the  abstract
Hardy-Sobolev spaces that we defined in the first part. Section \ref{section3} is devoted to the proof of the
interpolation results in Theorems  \ref{thm:hardy} and \ref{thm:hardyhomo}  using a ``Calder\'on-Zygmund''
decomposition well adapted to the spaces $HS^1_{(\B),ato}$.  For the interpolation of the abstract  Hardy-Sobolev
spaces in Theorem  \ref{thm:generalbis}, our method is  based on the new maximal inequality described in \cite{BB}.
Finally, the proof of Theorem \ref{HSL}  and the application to $\Delta^{1/2}$ are given in section \ref{section:App}.

\gb \paragraph{\textit{Acknowledgements:}} The two authors thank E.~Russ for his interest and the discussions they had
with him on this subject. The first author  would also like to thank F. Ricci for his hospitality during her visit to
the Scuola Normale di Pisa and for proposing her a part of this topic.
\section{Preliminaries}

\mb Throughout this paper we will denote by $\ind_{E}$ the characteristic function of
 a set $E$ and $E^{c}$ the complement of $E$. If $X$ is a metric space, $\Lip$ will be  the set of real Lipschitz functions on $X$ and $\Lip_{0}$ the set of real, compactly supported Lipschitz functions on $X$.  We denote by $Q(x, r)$ the open ball of
center $x\in X$ and radius $r>0$ and $\lambda Q$  denotes the ball co-centered with $Q$ and with radius $\lambda$ times that of $Q$. Finally, $C$ will be a constant that may change from an inequality to another and we will use $u\lesssim
v$ to say that there exists two constants $C$  such that $u\leq Cv$ and $u\simeq v$ to say that $u\lesssim v$ and $v\lesssim u$.

\mb In all this paper $M$ denotes a Riemannian manifold. We write $\mu$ for the Riemannian measure on $M$, $\nabla$ for the
Riemannian gradient, $|\cdot|$ for the length on the tangent space (forgetting the subscript $x$ for simplicity) and
$\|\cdot\|_{L^p}$ for the norm on $ L^p:=L^{p}(M,\mu)$, $1 \leq p\leq +\infty.$ \\
We will use the positive Laplace-Beltrami operator $\Delta$ defined by
$$ \forall f,g\in C^\infty_0(M), \qquad \langle \Delta f,g\rangle = \langle \nabla f ,\nabla g \rangle.$$

\subsection{The doubling property}

\begin{df} Let $M$ be a Riemannian manifold. One says that $M$ satisfies the (global) doubling property $(D)$ if there exists a
constant $C>0$, such that for all $x\in M,\, r>0 $ we have
\begin{equation*}\tag{$D$}
\mu(Q(x,2r))\leq C \mu(Q(x,r)).
\end{equation*}
\end{df}
\noindent Observe that if $M$ satisfies $(D)$ then
$$ \diam(M)<\infty\Leftrightarrow\,\mu(M)<\infty\,\textrm{ (see \cite{ambrosio1})}. $$
Therefore if $M$ is a complete non-compact Riemannian manifold satisfying $(D)$ then $\mu(M)=\infty$.

\begin{thm}[Maximal theorem]\label{MIT} (\cite{coifman2})
Let $M$ be a Riemannian manifold satisfying $(D)$. Denote by $ \mathcal{M}$ the uncentered Hardy-Littlewood maximal function
over open balls of $M$ defined by
 $$  \mathcal{M}f(x):=\underset{\genfrac{}{}{0pt}{}{Q \ \textrm{ball}}{x\in Q}} {\sup} \ |f|_{Q} $$
 where $\displaystyle f_{E}:=\aver{E}f d\mu:=\frac{1}{\mu(E)}\int_{E}f d\mu.$
Then for every  $p\in(1,\infty]$, $ \mathcal{M}$ is $L^p$ bounded and moreover of weak type $(1,1)$\footnote{ An operator $T$ is of weak type $(p,p)$ if there is $C>0$ such that for any $\alpha>0$, $\mu(\{x;\,|Tf(x)|>\alpha\})\leq \frac{C}{\alpha^p}\|f\|_p^p$.}.
\\
Consequently for $s\in(0,\infty)$, the operator $ \mathcal{M}_s$ defined by
$$  \mathcal{M}_sf(x):=\left[ \mathcal{M}(|f|^s)(x) \right]^{1/s} $$
is of weak type $(s,s)$ and $L^p$ bounded for all $p\in(s,\infty]$.
\end{thm}

\subsection{Poincar\'e inequality}
\begin{df}[Poincar\'{e} inequality on $M$] We say that a complete Riemannian manifold $M$ admits \textbf{a Poincar\'{e} inequality $(P_{q})$} for some $q\in[1,\infty)$ if there exists a constant $C>0$ such that, for every function $f\in \Lip_{0}(M)$\footnote{compaclty supported Lipshitz function defined on $M$.} and every ball $Q$ of $M$ of radius $r>0$, we have
\begin{equation*}\tag{$P_{q}$}
\left(\aver{Q}|f-f_{Q}|^{q} d\mu\right)^{1/q} \leq C r \left(\aver{Q}|\nabla f|^{q}d\mu\right)^{1/q}.
\end{equation*}
\end{df}
\begin{rem} By density of $C_{0}^{\infty}(M)$ in $Lip_0(M)$, we can replace $\Lip_{0}(M)$ by $C_{0}^{\infty}(M)$.
\end{rem}
Let us recall some known facts about Poincar\'{e} inequalities with varying $q$.
 \\
It is known that $(P_{q})$ implies $(P_{p})$ when $p\geq q$ (see \cite{hajlasz4}). Thus if the set of $q$ such that
$(P_{q})$ holds is not empty, then it is an interval unbounded on the right. A recent result of S. Keith and X. Zhong
(see \cite{KZ}) asserts that this interval is open in $[1,+\infty[$~:

\begin{thm}\label{kz} Let $(X,d,\mu)$ be a complete metric-measure space with $\mu$ doubling
and admitting a Poincar\'{e} inequality $(P_{q})$, for  some $1< q<\infty$.
Then there exists $\epsilon >0$ such that $(X,d,\mu)$ admits
$(P_{p})$ for every $p>q-\epsilon$.
\end{thm}

\mb A consequence of Poincar\'e inequality:

\begin{prop} \label{prop:moy} Assume that $M$ satisfies $(D)$ and  admits a Poincar\'e inequality $(P_p)$  for some $p\in[1,\infty)$. Then there is a constant $c=c(p)$ such that for all balls $Q$ (of radius $r_Q$) and all functions $f\in C^\infty_0(Q)$
$$ \left| \frac{1}{\mu(Q)}\int_Q f d\mu \right| \leq c r_Q \left( \frac{1}{\mu(Q)}\int_Q |\nabla f|^{p} d\mu \right)^{1/p}.$$
\end{prop}

\mb This result is well-known. However for an easy reference and for the sake of completeness, we remember the proof
based on the self-improvement of Poincar\'e inequality. We refer the reader to Theorem 5.3.3 of \cite{SC} for an
initial proof (the proof there applies also for $p=1$).

\mb\dem We first prove that for all $x\in Q, y\in 3Q\setminus 2Q$ \be{prop:moy2} \left|f(x)-f(y)\right| \lesssim r_Q
 \mathcal{M}_{p-\epsilon}(|\nabla f|)(x). \ee  Using Hardy-Littlewood Theorem, we have
$$ f(x) = \lim_{\epsilon \to 0} f_{Q(x,\epsilon)}.$$
With the balls $Q_i:=Q(x,2^ir_Q)$, we also have
$$ \left| f(x)-f_{Q_1} \right| \leq \sum_{i\leq 0} \left| f_{Q_i} - f_{Q_{i+1}}\right|.$$
Thanks to Theorem \ref{kz}), the Poincar\'e inequality $(P_p)$ self improves to $(P_{p-\epsilon})$ for a certain $\epsilon>0$. Using this Poincar\'e inequality and the doubling property one obtains
\begin{align*}
 \left| f(x)-f_{Q_4} \right| & \leq \sum_{i=-\infty}^{3} \left| f_{Q_i} - f_{Q_{i+1}}\right| \\
 & \lesssim \sum_{i\leq 3} \frac{1}{\mu(Q_i)}\int_{Q_i} \left| f - f_{Q_i} \right| d\mu \\
 & \lesssim \sum_{i\leq 3} r_{Q_i} \left(\frac{1}{\mu(Q_i)}\int_{Q_i} \left| \nabla f \right|^{p-\epsilon} d\mu\right)^{\frac{1}{p-\epsilon}} \\
 & \lesssim \sum_{i\leq 3} 2^{-i}r_Q  \mathcal{M}_{p-\epsilon}(|\nabla f|)(x) \\
 & \lesssim r_Q  \mathcal{M}_{p-\epsilon}(|\nabla f|)(x).
\end{align*}
Similarly we have with $\widetilde{Q}_i:=Q(y,2^ir_Q)$
\begin{align*}
  \left| f(y)-f_{\widetilde{Q}_3} \right| \lesssim r_Q  \mathcal{M}_{p-\epsilon}(|\nabla f|)(y).
\end{align*}
However since $y\in 3Q\setminus 2Q$ and $f$ is supported in $Q$, we have $ \mathcal{M}_{p-\epsilon}(|\nabla f|)(y)\lesssim
 \mathcal{M}_{p-\epsilon}(|\nabla f|)(x)$. Then we just have to control the difference of means. The Poincar\'e inequality
$(P_p)$ and $\widetilde{Q}_3\subset Q_4$ yield
$$ \left|f_{\widetilde{Q}_3} - f_{Q_4} \right|\lesssim \frac{1}{\mu(Q_4)}\int_{Q_4} \left|f-f_{Q_4}\right| d\mu \lesssim r_Q   \mathcal{M}_{p-\epsilon}(|\nabla f|)(x).$$
Thus  we proved (\ref{prop:moy2}). Then using the fact that $f(y)=0$ due to the support of $f$, we obtain~
$$ \left| \frac{1}{\mu(Q)}\int_Q f d\mu \right| \leq \frac{1}{\mu(Q)} \int_Q \left|f(x)-f(y)\right| d\mu(x) \lesssim r_Q\mu(Q)^{-1/p} \left\|  \mathcal{M}_{p-\epsilon}(|\nabla f |) \right\|_{p}.$$
Finally the $L^p$-boundedness of $ \mathcal{M}_{p-\epsilon}$ concludes the proof. \findem

\subsection{The $K$-method of real interpolation}
\label{subsec:rappel}

The reader can refer to \cite{bennett}, \cite{bergh} for details on the development of this theory. Here we only recall the essentials to be used in the sequel.

\mb Let $A_{0}$, $A_{1}$ be  two normed vector spaces embedded in a topological Hausdorff vector space $V$. For each  $a\in A_{0}+A_{1}$ and $t>0$, we define the $K$-functional of real interpolation by
$$
K(a,t,A_{0},A_{1})=\displaystyle \inf_{a
=a_{0}+a_{1}}(\| a_{0}\|_{A_{0}}+t\|
a_{1}\|_{A_{1}}).
$$

\gb For $0<\theta< 1$, $1\leq q\leq \infty$, we denote by $(A_{0},A_{1})_{\theta,q}$ the real interpolation space between $A_{0}$ and $A_{1}$ defined as
\begin{displaymath}
    (A_{0},A_{1})_{\theta,q}=\left\lbrace a \in A_{0}+A_{1}:\|a\|_{\theta,q}=\left(\int_{0}^{\infty}(t^{-\theta}K(a,t,A_{0},A_{1}))^{q}\,\frac{dt}{t}\right)^{\frac{1}{q}}<\infty\right\rbrace.
\end{displaymath}
It is an exact interpolation space of exponent $\theta$ between $A_{0}$ and $A_{1}$ (see \cite{bergh}, Chapter II).
\begin{df}
Let $f$  be a measurable function on a measure space $(X,\mu)$. The decreasing rearrangement of $f$ is the function $f^{*}$ defined for every $t\geq 0$ by
$$
f^{*}(t)=\inf \left\lbrace\lambda :\, \mu (\left\lbrace x:\,|f(x)|>\lambda\right\rbrace)\leq
t\right\rbrace.
$$
The maximal decreasing rearrangement of
$f$ is the function $f^{**}$ defined for every $t>0$ by
$$
f^{**}(t)=\frac{1}{t}\int_{0}^{t}f^{*}(s) ds.
$$
\end{df}
From the properties of $f^{**}$ we mention:
\begin{itemize}
\item[1.] $(f+g)^{**}\leq f^{**}+g^{**}$.
\item[2.] $(\mathcal{M}f)^{*}\sim f^{**}$.
\item[3.] $\mu(\left\lbrace x;\,|f(x)|>f^{*}(t)\right\rbrace )\leq t$.
\item[4.] $\forall 1< p\leq\infty$, $\|f^{**}\|_{p}\sim\|f\|_{p}$.
\end{itemize}

\mb We exactly know the functional $K$ for Lebesgue spaces~:
\begin{prop} Take $0<p_0<p_1\leq \infty$. We have~:
$$K(f,t,L^{p_0},L^{p_1}) \simeq \left(\int_0^{t^{\alpha}} \left[f^{*}(s)\right]^{p_0} ds \right)^{1/p_0} + t \left(\int_{t^{\alpha}}^\infty \left[f^{*}(s)\right]^{p_1} ds \right)^{1/p_1},$$
where $\frac{1}{\alpha}=\frac{1}{p_0}-\frac{1}{p_1}$.
\end{prop}

\mb From now on, we always assume that the Riemannian manifold satisfies the doubling property $(D)$.

\subsection{Maximal inequalities for dual Sobolev spaces.} \label{subsec:maximal}

\mb First, we begin recalling the ``duality-properties'' of the Sobolev spaces.

\begin{df} For $p\in[1,\infty]$ and $O$ an open set of $M$, we define $W^{1,p}(O)$ as following
$$ W^{1,p}(O) := \overline{C^\infty_0(O)}^{\|\, .\, \|_{W^{1,p}(O)}}, \qquad \textrm{with} \qquad \|f\|_{W^{1,p}(O)}:=\left\||f|+|\nabla f| \right\|_{L^p(O)}.$$
Then we denote $W^{-1,p'}(O)$ the dual space of $W^{1,p}(O)$ defined as the set of distributions $f\in \mathcal{D}'(M)$ such that
$$
\|f\|_{W^{-1,p'}(O)}=\sup_{g\in C^{\infty}_{0}(M)} \frac{|\langle f,g\rangle|}{\|g\|_{W^{1,p}(O)}}.
$$
\end{df}

\begin{prop} \label{dualS} Let $p\in[1,\infty)$. Then for all open set $O$ of $M$, we have
\begin{align*}
 \|f\|_{W^{-1,p'}(O)} & \simeq \inf_{f=\phi-div(\psi)} \left\|\phi\right\|_{L^{p'}(O)} + \left\| \psi\right\|_{L^{p'}(O)} \\
 & \simeq \inf_{f=\phi-div(\psi)} \left\||\phi|+|\psi|\right\|_{L^{p'}(O)}.
\end{align*}
Here we take the infimum over all the decompositions $f=\phi-div(\psi)$ on $M$ with $\phi\in L^{p'}(O)$ and $\psi\in
\D'(O,\R^n)$ such that $div(\psi)\in L^{p'}(O)$.
\end{prop}

\mb The proof is left to the reader (it is essentially written in \cite{ART}, Proposition 33).

\mb We now introduce the following maximal operators~:

\begin{df} Let $s>0$. According to the standard maximal ``Hardy-Littlewood'' operator
$  \mathcal{M}_s$,
we define two ``Sobolev versions''~:
$$ M_{S,s}(f)(x):=\sup_{\genfrac{}{}{0pt}{}{Q \textrm{ball}}{x\in Q}} \ \frac{1}{\mu(Q)^{1/s}} \left\|f \right\|_{W^{-1,s}(Q)}$$
and
$$ M_{S,*,s}(f)(x):= \inf_{f = \phi - div(\psi)} \  \mathcal{M}_s\left(|\phi|+|\psi|\right)(x).$$
\end{df}

\mb The following assumption is taken from  \cite{BB}:

\begin{ass} Take two exponents $1\leq\p_0\leq \p_1<\infty$. We call $(H_{\p_0,\p_1})$ the following assumption~:
 \begin{equation}
\|f\|_{W^{-1,\p_1}} \lesssim \|M_{S,*,\p_0}(f)\|_{L^{\p_1}}. \tag{$H_{\p_0,\p_1}$} \label{Hp}
\end{equation}
\end{ass}

\begin{df} \label{ensembleI} For $M$ a Riemannian manifold, we denote by ${\mathcal I}_M$ the following set

$$ {\mathcal I}_M:=\left\{ (\p_0,\p_1)\in(1,\infty)^2,\ \p_0\leq \p_1,\ (\textrm{\ref{Hp})} {\textrm{ holds}} \right\}.$$
\end{df}

\mb We refer to \cite{BB} for the study of these maximal operators and the previous assumption.

\begin{prop} \label{prop:importante}
For $p\in[1,\infty)$, $M_{S,p}$ and $M_{S,*,p}$ are of ``weak type $(p,p)$''. That is \be{eq1} \forall f\in
W^{-1,p},\qquad  \left\| M_{S,p}(f) \right\|_{L^{p,\infty}} \leq \left\| M_{S,*,p}(f) \right\|_{L^{p,\infty}} \lesssim
\|f\|_{W^{-1,p}}. \ee
\end{prop}

\begin{df} We use the operator $L:=(I+\Delta)$ defined with the positive Laplace-Beltrami operator. We recall that the two operators $\Delta$ and $L$ are self-adjoint. \\
 According to \cite{AC}, we say that for $p\in(1,\infty)$ we have the non-homogeneous property (\ref{Rp}) if
 \begin{equation} \label{Rp}
 \|f\|_{W^{1,p}} \lesssim \left\| L^{1/2}(f) \right\|_{L^p} \tag{$nhR_p$}
 \end{equation}
 for all $f\in C_0^{\infty}(M)$. This is equivalent to the $L^p$ boundedness of the local Riesz transform $\nabla(I+\Delta)^{-1/2}$.\
We have the non-homogeneous reverse property (\ref{RRp}) if
 \begin{equation} \label{RRp}
   \left\| L^{1/2}(f) \right\|_{L^p} \lesssim \|f\|_{W^{1,p}} \tag{$nhRR_p$}
 \end{equation}
  for all $f\in C_0^{\infty}(M)$.
\end{df}

\begin{df} \label{df:od} Let $p,q\in[1,\infty)$. We say that the collection $(T_t)_{t>0}=(e^{-t\Delta})_{t>0}$ or $(T_t)_{t>0}=(\sqrt{t}\nabla e^{-t\Delta})_{t>0}$  satisfy $(L^{p}-L^{q})$-``off-diagonal'' estimates, if there exists $\gamma$ such that for all balls $Q$ of radius $r_Q$, every function $f$ supported in $Q$ and all index $j\geq 0$
$$ \left( \frac{1}{\mu(2^jQ)}\int_{S_j(Q)} \left| T_{r_Q^2}(f) \right|^{q} d\mu \right)^{1/q} \lesssim  e^{-\gamma 4^j} \left( \frac{1}{\mu(Q)}\int_Q \left| f \right|^{p} d\mu \right)^{1/p}.$$
We used $S_j(Q)$ for the dyadic corona around the ball
$$S_j(Q):=\left\{y, 2^{j}\leq 1+\frac{d(y,Q)}{r_Q} <2^{j+1}\right\}.$$
These ``off-diagonal'' estimates are closely related to ``Gaffney estimates'' of the semigroup.
\end{df}

\mb We now come to the main result of \cite{BB}.

\begin{thm} \label{poincare} Let $1<s<r'<\sigma$. Assume that the Riemannian manifold $M$ satisfies $(nhRR_{r})$ and $(nhR_{s'})$. Moreover, assume that the semigroup $(e^{-t\Delta})_{t>0}$ satisfies $(L^{\sigma'}-L^{s'})$-``off-diagonal'' estimates and that the collection $(\sqrt{t}\nabla e^{-t\Delta})_{t>0}$ satisfies $(L^{s'}-L^{s'})$-``off-diagonal'' estimates.
Then there is a constant $c=c(s,r,\sigma)$ such that
\be{eq3bisse} \forall f\in W^{-1,r'}, \qquad \|f\|_{W^{-1,r'}} \lesssim \left\| M_{S,*,s}(f) \right \|_{L^{r'}}. \ee
Therefore (\ref{Hp}) is satisfied for all exponents $\p_0,\p_1$ satisfying $ \p_0\geq s$ and $\p_1 =r'$.
\end{thm}

\begin{cor} \label{corRn} In the Euclidean case $M=\R^n$, for all $\p_0,\p_1\in(1,\infty)$, the assumption (\ref{Hp}) holds. More generally, on any Riemannian manifold satisfying $(D)$ and $(P_1)$, (\ref{Hp}) holds for all $\p_0,\p_1\in(1,\infty)$.
\end{cor}

\mb After all these preliminaries, we now define our Hardy-Sobolev spaces via atomic decomposition.

\section{Abstract Hardy-Sobolev spaces.}
\label{sec:method1}
 We begin this section defining ``abstract atomic''
Hardy-Sobolev spaces, then we study in more detail a particular case of these spaces.

\subsection{New Hardy-Sobolev spaces.} \label{subsec:method1}

We follow ideas of \cite{BJ} and propose an ``atomic'' definition of abstract Hardy-Sobolev spaces. We refer the reader to \cite{BJ} for an explanation of this choice~: the ``atoms'' are defined as the image of localized functions by an operator $B_Q$, playing the role of the ``oscillation operator'' associated to a ball $Q$.

\gb Let us fix $\B\in]1,\infty]$ and take  $\mathbb{B}:=(B_Q)_{Q\in \mathcal{Q}}$ a collection of $W^{1,\B}$-bounded linear operators, indexed by $\mathcal{Q}$ the collection of all open balls $Q$ of the manifold $M$. We assume that these operators $B_Q$ are uniformly bounded on $W^{1,\B}$ : there exists a
constant $0<A'<\infty$
such that~ \be{operh} \forall f\in W^{1,\B} ,\ \forall  Q \textrm{ ball}, \qquad \|B_Q(f) \|_{W^{1,\B}} \leq A'\|f\|_{W^{1,\B}}. \ee

\mb We define the {\it Sobolev-atoms}  using the collection $\mathbb{B}$~:

\begin{df} A function $m\in L^{1}_{loc}$ is called an atom associated to a ball $Q$ if there exists a real function $f_Q$ compactly supported in the ball $Q$ such that
$$m=B_Q(f_Q),$$
with
$$ \left\|f_Q\right\|_{W^{1,\B}} \leq \mu(Q)^{-1/\B'}.$$
\end{df}

\noindent The functions $f_Q$ in this definition are normalized in
$W^{1,1}$. It is easy to check that
$$ \|f_Q\|_{W^{1,1}} \lesssim 1 .$$
Now we can define our abstract atomic Hardy-Sobolev spaces~:

\begin{df} \label{df1} A measurable function $h$ belongs to the atomic
Hardy-Sobolev space $HW^1_{ato}$ if there exists a decomposition~
$$h=\sum_{i\in {\mathbb N}} \lambda_i m_i  \qquad \mu-a.e, $$
where for all $i$, $m_i$ is an atom and $(\lambda_{i})_i$
are real numbers satisfying
$$\sum_{i\in {\mathbb N}} |\lambda_i| <\infty. $$
We equip $HW^1_{ato}$ with the norm~:
$$\|h\|_{HW^1_{ato}}:= \inf_{h=\sum_{i\in {\mathbb N}} \lambda_i m_i} \sum_{i} |\lambda_i|.$$
Similarly we define our ``finite'' Hardy-Sobolev space $HW^1_{F,ato}$ as the set of functions which admit finite atomic decompositions.
\end{df}

\begin{rem} We refer the reader to \cite{BJ,BCRASS,B2} for details concerning the use of ``finite atomic Hardy space'' instead of the whole atomic Hardy space. The use of this last one brings technical problems (we do not know how to solve them) that are not important and are twisted by the use of the atomic Hardy space.
\end{rem}

\mb Our goal is to interpolate the Hardy-Sobolev spaces with Sobolev spaces. First, we  describe a useful criterion to prove the boundedness of an operator from the Hardy-Sobolev space $HW^{1}_{F,ato}$ into $L^1$.

\begin{prop} \label{prop:contHW} Let $M$ be a Riemannian manifold satisfying the doubling property. Let $T$ be a linear operator  bounded from $W^{1,\B}$ to $L^\B$  for some $\B \in (1,\infty)$  and
satisfying some ``off-diagonal'' Sobolev estimates: for all ball $Q$ and all function $f$ compactly supported in $Q$
\be{offdiagonal} \forall j\geq 2 \quad \left(\frac{1}{\mu(2^{j+1}Q)} \int_{S_j(Q)} \left| T(B_Q(f))\right|^{\B} d\mu
\right)^{1/\B} \leq \alpha_j(Q) \frac{1}{\mu(Q)^{1/\B}} \|f\|_{W^{1,\B}(Q)}, \ee with coefficients $\alpha_j$
satisfying \be{assum:10} \Lambda:=\sup_{Q \textrm{ ball}} \  \sum_{j\geq 2} \frac{\mu(2^{j+1}Q)}{\mu(Q)} \alpha_j(Q)
<\infty. \ee Then $T$ is continuous from $HW^1_{F,ato}$ to $L^1$.
\end{prop}

\mb The proof is left to the reader, it is written in \cite{BJ} and \cite{B2} in the context of Lebesgue spaces. It is
the same in our context of Sobolev spaces.

\subsection{The study of a particular Hardy-Sobolev space.} \label{sub:par}

In this subsection, we present in more detail the study of a particular Hardy-Sobolev space.

\mb In the study of Hardy spaces (see \cite{BJ}), we have seen that our abstract Hardy space corresponds to the ``classical'' Hardy space (the one defined by R. Coifman and G. Weiss in \cite{coifman1}), when we choose our operator $B_Q$ as the exact oscillation operator. Here we want to study the Hardy-Sobolev space defined with a regular version of this particular collection ${\mathbb B}$. For all ball $Q$, let $\phi_Q$ be a function supported in $Q$ and satisfying~
$$ \|\phi_Q \|_{L^\infty}\lesssim 1, \quad \|\,| \nabla \phi_Q |\,\|_{L^\infty} \lesssim r_Q^{-1} \quad \textrm{and} \quad \int \phi_Q d\mu =\mu(Q).$$
We define our operator
$$ A_Q(f):=\left(\frac{1}{\mu(Q)}\int_Q f d\mu \right) \phi_Q \qquad B_Q(f)=f-A_Q(f).$$
In all this subsection, the Hardy-Sobolev spaces are constructed with this particular choice of operators. According to this collection, we construct our Hardy-Sobolev space $HW^{1}_{(\B),ato}$ and $HW^{1}_{F,(\B),ato}$.

\begin{rem} In the previous subsection, we did not study the dependence of the Hardy-Sobolev space with respect to the exponent $\B$, so we omitted it in the notation. In this subsection, we will study the role of $\B$ in a particular case (see Theorem \ref{thm:comp}). That is why we put the exponent in the notation.
\end{rem}

\mb We have to check the first assumption (\ref{operh}). Thanks to Proposition \ref{prop:moy}, it is easy to check that if a Poincar\'e inequality $(P_{\B})$ is satisfied then (\ref{operh}) holds.

\mb Moreover, with the normalization of functions $\phi_Q$, each atom $m$ associated to a ball $Q$ verifies
$$\int_Q m d\mu = 0.$$
From this observation, we can set a definition of particular Hardy-Sobolev spaces.

\begin{df} For $\B\in(1,\infty]$, we say that a function $m$ is a non-homogeneous $(1,\B)$-atom associated to a ball $Q$, if
\begin{itemize}
\item[1.] $m$ is supported in the ball $Q$,
\item[2.] $\|m\|_{W^{1,\B}}\leq\mu(Q)^{-\frac{1}{\B'}}$,
\item[3.] $\int m d\mu=0$.
\end{itemize}
We define the Hardy-Sobolev space $HS^{1}_{(\B),ato}$ as follows: $f\in HS^{1}_{(\B),ato}$ if there exists $(b_{i})_{i}$ a family of $(1,\B)$-atoms such that $f=\sum_{i}\lambda_{i}b_{i}$ with $\sum_{i}|\lambda_{i}|<\infty$.
We equip this space with the norm
$$
\|f\|_{HS^{1}_{(\B),ato}}=\inf_{(\lambda_{i})_{i}} \sum_{i}|\lambda_{i}|.
$$
Similarly to Definition \ref{df1}, we define ``finite'' atomic space
$HS^{1}_{F,(\B),ato}$.
\end{df}

\mb From Proposition \ref{prop:moy} and the previous discussion, we have this first proposition.

\begin{prop} Assume that a Poincar\'e inequality $(P_\B)$ holds. Then the concept of $(1,\B)$-atoms exactly corresponds to the concept of atoms, defined with our operators $B_Q$. Thus the different atomic Hardy-Sobolev spaces are equal:
$$ HS^{1}_{F,(\B),ato}=HW^{1}_{F,(\B),ato} \qquad HS^{1}_{(\B),ato}=HW^{1}_{(\B),ato}.$$
\end{prop}

\begin{rem} Note that every $\B_2$ atom is an $\B_1$ atom for $1<\B_1\leq \B_2\leq\infty$ and therefore $HS_{(\B_2),ato}^{1}\subset HS_{(\B_1),ato}^{1}$ with $\|f\|_{HS^{1}_{(\B_1),ato}}\leq \|f\|_{HS_{(\B_2),ato}^{1}}$.
\end{rem}

\begin{prop} $HS_{(\B),ato}^1$ is a Banach space for $\B\in(1,\infty]$.
\end{prop}
\dem Consider a sequence $(h_{k})_k$ in $HS_{(\B),ato}^{1}$ such that $\sum_k\|h_k\|_{HS_{(\B),ato}^1}<\infty$. It suffices to prove that $\sum_k h_k$ converges in $HS_{(\B),ato}^1$. For this, for every $k$ take the following atomic decomposition $h_k=\sum_i \lambda_{k,i}b_{k,i}$ with $\sum_i|\lambda_{k,i}|\leq \|h_k\|_{HS_{(\B),ato}^{1}}+\frac{1}{2^k}$. Then $h=\sum_k\sum_i\lambda_{k,i}b_{k,i}\in W_1^1$ (absolutely convergence) with $\sum_{k}\sum_i|\lambda_{k,i}|\leq \sum_{k}\|h_k\|_{HS_{(\B),ato}^{1}}+\sum_{k}\frac{1}{2^k}<\infty$. Hence $h\in HS_{(\B),ato}^1$ and the proof is complete. \findem

\begin{prop}\label{DS} For $\B\in(1,\infty]$, the finite space $HS_{F,(\B),ato}^1$ is dense in $HS_{(\B),ato}^1$.
\end{prop}
We recall here the definition of a Coifman-Weiss atom of $H_{CW}^1:=H^1_{CW}(M)$  the Hardy
space of Coifman-Weiss (see \cite{coifman1}).

\begin{df} \label{atCW}For $\B\in(1,\infty]$, we say that a function $m$ is a $\B$-atom associated to a ball $Q$, if
\begin{itemize}
\item[1.] $m$ is supported in the ball $Q$,
\item[2.] $\|m\|_{L^{\B}}\leq\mu(Q)^{-\frac{1}{\B'}}$,
\item[3.] $\int m d\mu=0$.
\end{itemize}
\end{df}
In the literature, we found  definitions of classical Hardy-Sobolev spaces in the Euclidean case as the set of $f\in
H^{1}_{CW}$ such that $\nabla f \in H^{1}_{CW}$ or $\Delta^{1/2}f \in H^{1}_{CW}$. Thanks to the $H_{CW}^1$ boundedness
of the Riesz transform  in $\mathbb{R}^{n}$, these two spaces are equal.

We hope to have a complete picture and comparison of all these definitions of Hardy-Sobolev spaces on Riemannian
manifolds in a forthcoming paper.
\begin{df} The classical Hardy-Sobolev space $HS^{1}(M)$ is defined as (see \cite{cho},
for the Euclidean case)
$$HS^{1}=\left\{ f\in H_{CW}^{1};\, \nabla f \in H^1_{CW}
 \right\}$$
 where $\nabla f$ is the distributional gradient of $f$.
 \end{df}

\begin{prop} The space $HS^{1}$ is a Banach space. \end{prop}
\dem
Let $(f_n)_n$ be a Cauchy sequence in $HS^{1}$. Then $(f_{n})$ and $(\nabla f_{n})_{n}$ are Cauchy sequences in $H_{CW}^{1}$ and therefore converge to $f\in H_{CW}^{1}$  and $g\in H_{CW}^{1}$. Since $f_{n}\rightarrow f\; \mu-a.e$ it comes that $\nabla f_{n}\rightarrow \nabla f$ in the distributional sense. The uniqueness of the limit shows that $g=\nabla f$ and finishes the proof.
\findem

\begin{prop} We have $HS_{(\B),ato}^{1} \subset HS^{1} \subset W^{1,1}$ for every $\B>1$.
\end{prop}

\mb Unfortunately, it is not clear when $HS^{1}\subset HS_{(\B),ato}^{1}  $. However for the point of view of
interpolation, the study of $HS_{(\B),ato}^1$ implies results for $HS^{1}$. For an exponent $p\in(1,\infty]$ and
$\theta\in(0,1)$ if
$$\left(HS_{(\B),ato}^1,W^{1,p}\right)_{\theta,p_\theta} = W^{1,p_\theta}$$ with
$\frac{1}{p_\theta}=(1-\theta) + \frac{\theta}{p}$ then we know that $$\left(HS^1,W^{1,p}\right)_{\theta,p_\theta} =
W^{1,p_\theta}.$$
This follows from the fact that $HS_{(\B),ato}^{1} \subset HS^{1} \subset W^{1,1}$ and $\|f\|_{HS^1}\leq 2 \|f\|_{HS_{(\B),ato}^1}$.

\gb  We know (see \cite{coifman1}) that the Hardy space $H^1_{CW}$
admits an atomic decomposition and is also equal to the corresponding atomic Hardy space (for any exponent $\B$ used in
the definition of $\B$-atoms). In our case the atomic Hardy-Sobolev spaces are all contained in the classical one
$HS^1$ but for the moment we are not able to show if they are equal or not. We believe that this is not true without additional hypotheses on the geometry of the manifold. However, under Poincar\'e inequality we
will compare different atomic Hardy-Sobolev spaces in Theorem \ref{thm:comp}.

\mb Before we prove this theorem, we need the following Lemma.
\begin{lem} \label{coif}(see Lemma 3.9 in \cite{coifman1}) Assume that $M$ satisfies $(D)$.
\begin{itemize}
\item[1.] Let $$  \mathcal{M}_{c}f(x):=\sup_{r>0}\frac{1}{\mu(Q(x,r))}\int_{Q(x,r)}|f|d\mu$$ be the centered maximal function of $f$. Observe that if $x\in Q(y,r)$ then $Q(y,r)\subset Q(x,2r)$. It follows that
 $$
  \mathcal{M}_{c}f\leq  \mathcal{M}f\leq C \mathcal{M}_{c}f
 $$
 where $C$ only depends on the constant of the doubling property.
 \item[2.] Let $f$ be an $L^{1}$ function supported in $Q_0=Q(x_0,r_0)$. Then there is $C_1$ depending on the doubling constant such that
 $$\Omega_{\alpha}=\left\lbrace x\in M;  \mathcal{M}(f)(x)>\alpha\right\rbrace \subset Q(x_0,2r_0) $$
 whenever $\alpha>C_{1}\avert{Q_{0}}|f|d\mu$.
 \end{itemize}
\end{lem}
{\bf Proof of Theorem \ref{thm:comp}~:} \\
The proof is inspired by that of R. Coifman and G. Weiss (\cite{coifman1}) for classical Hardy spaces on a space of
homogeneous type. We prove that every $(1,\B)$-atom is a sum of $(1,\infty)$-atoms. We use an adapted Calder\'{o}n-Zygmund
decomposition for Sobolev functions (proved later for convenience in Subsection \ref{sub:intpar}) and proceed as their
proof. However, the presence of the gradient create some problems. \\ Since we  know from \cite{KZ} that Poincar\'e inequality $(P_q)$ self-improves in $(P_{q-\epsilon})$ for some
$\epsilon>0$,  let us denote $\kappa:=\| \mathcal{M}_{q-\epsilon}\|_{L^\beta\to L^\beta}$. \\
Let $a$ be a $(1,\B)$ atom supported in a ball $Q_0$. Set $b=\mu(Q_{0})$a. \\
We claim that for $K,\alpha>0$ large enough parameters and numerical constants $C$ and $N$, there exists a collection of balls $(Q_{j_l})$, $j_l\in
\mathbb{N}^{l}$ for $l=0,1,...,$ such that for every $n\geq 1$
 \begin{equation}\label{bex}
 b=CNK\alpha \sum_{l=0}^{n-1}(K\alpha)^{l}\sum_{j_{l}\in \mathbb{N}^{l}}\mu(Q_{j_l}) a_{j_l}+\sum_{j_n\in \mathbb{N}^{n}}h_{j_n}
 \end{equation}
 and
 \begin{itemize}
 \item[(a)]
$ a_{j_l}$ is an $\infty$ atom supported in $Q_{j_l}$, $l=0,1,...n-1$;
\item[(b)]$\bigcup_{j_{n}\in
\mathbb{N}^{n}}Q_{j_{n}}\subset \Omega_{n}:=\left\lbrace x;\, \mathcal{M}_{q}\left(|b|+ \mathcal{M}_{q-\epsilon}^n(|\nabla b|)
\right)(x)>K\frac{\alpha^{n}}{2}\right\rbrace$;
 \item[(c)] $\sum_{j_n}\ind_{Q_{j_n}}\leq N^{n}$;
 \item[(d)] $\supp h_{j_n}\subset Q_{j_n}$, $ \,\int h_{j_n}d\mu=0$;
 \item[(e)] $|h_{j_n}(x)|\leq |b(x)|+ 2C'(K\alpha)^{n}\ind_{Q_{j_n}}(x)$;
 \item[(f)] $|\nabla h_{j_n}(x)|\leq (C''+2) \left[K^{n-1}  \mathcal{M}_{q-\epsilon}^n(|\nabla b|)(x)+ (K\alpha)^{n}\ind_{Q_{j_n}}(x)\right]$;
 \item[(g)] $\left(\avert{Q_{j_n}} (|h_{j_n}|^{q}+|\nabla h_{j_n}|^{q})d\mu \right)^{\frac{1}{q}}\lesssim (K\alpha)^{n}$.
 \end{itemize}
 The constants $\alpha, K$  are sufficiently large and $\alpha,K, \, N$ depend only on $\B, \,q$ and the doubling constant.
 We write $ \mathcal{M}_{q-\epsilon}^n$ for the composed operator $ \mathcal{M}_{q-\epsilon} \circ  \mathcal{M}_{q-\epsilon} \circ ...$.  \\
 Let us first see how from theses properties we can write
 $$a=\sum_{j}\alpha_{j}a_{j}$$
 where for every $j$, $a_{j}$ is an $\infty$-atom.
 We have
\begin{equation}\label{salpha}
\frac{NK\alpha}{\mu(Q_0)}\sum_{n=0}^{\infty}(K\alpha)^{n}\sum_{j_n\in \mathbb{N}^{n}}\mu(Q_{j_n})\leq C
 \end{equation}
 where $C$ is independent of $a$. Indeed, it follows from (b),\,(c) and the weak type $(1,1)$ of $ \mathcal{M}$ that
 \begin{equation*}\tag{$3'$}
 \sum_{j_n}\mu(Q_{j_n})\leq CN^{n}\mu(\bigcup_{j_n}Q_{j_n})\leq CN^n\mu(\Omega_n)\leq C\max(1,\kappa^n)N^n\left(\frac{2}{\alpha^{n}}\right)^{\B}\|b\|_{W^{1,\B}}^{\B}.
 \end{equation*}
 Therefore
 $$
 \sum_{n=0}^{\infty}(K\alpha)^{n}\sum_{j_n\in \mathbb{N}^{n}}\mu(Q_{j_n})\leq C2^\B\sum_{n=0}^{\infty}(\max(1,\kappa) NK\alpha^{1-\B})^{n}\|b\|_{W^{1,\B}}^{\B}.
 $$
Since $\|b\|_{W^{1,\B}}^{\B}\leq C\mu(Q_{0})$ we deduce (\ref{salpha}) with $C$ depending on $\B,\,q$ and $\alpha, K$ but
not on $a$. We choose $\alpha>>K$ such that $\frac{N\max(1,\kappa)K}{\alpha^{\B-1}}<1$.  \\
 From (g), we have
 $$
 \int (|h_{j_n}|+|\nabla h_{j_n}|)d\mu\leq C\mu(Q_{j_n})(K\alpha)^{n}.
 $$
 Therefore, if we note by $H_n=\sum_{j_n\in \mathbb{N}^{n}}h_{j_n}$ we have
 \begin{align*}\int(|H_n|+|\nabla H_n|)d\mu &\leq \sum_{j_n\in \mathbb{N}^{n}}\int(|h_{j_n}|+|\nabla h_{j_n}|)d\mu
 \\
 &\leq C(K\alpha)^{n}\sum_{j_{n}}\mu(Q_{j_n})
 \\
 &\leq C \left(N\max(1,\kappa) K\alpha^{1-\B} \right)^{n}\|b\|_{W^{1,\B}}^{\B}.
 \end{align*}
 We used the bounded overlap property of the $(Q_{j_n})_{n}$ and the above for $\mu(Q_{j_n})$. This shows that the first series on the right-hand side of (\ref{bex}) converges to $b$ in $W^{1,1}$. \\
 It remains now to prove that these properties are valid for every $n\in \mathbb{N}^{*}$. We begin proving the case $n=1$. Let
 $$\widetilde{\Omega_1}:= \left\lbrace x; \mathcal{M}(|b|^q+|\nabla b |^q)(x)>(K\alpha)^q\right\rbrace \subset \Omega_1.$$
 Lemma \ref{coif} shows that  $\widetilde{\Omega_{1}}\subset 2Q_{0}$ provided $K\alpha>C_1$. Moreover
$\widetilde{\Omega_{1}}$ is a bounded open set with $\mu(\widetilde{\Omega_{1}})\leq
\frac{C}{\alpha^{\B}}\|b\|_{W^{1,\B}}^{\B}\leq C \alpha^{-\B}\mu(Q_0)$. This allows us to apply the Whitney covering
theorem to $\widetilde{\Omega_{1}}$ and consider the Calder\'on-Zygmund decomposition of Proposition \ref{CZ} -- in section 3-- for $b$
with $p=\B$.
 We obtain
 \begin{equation}\label{bhg} b=\sum_{j}h_{j}+g_0
 \end{equation}
  with $h_{j}$, $g_{0}$ satisfying the properties of Proposition \ref{CZ}.
 We have
 \begin{align*} \sum_{j}\int_{Q_{j}}(|h_{j}|+|\nabla h_j|)d\mu&\leq C\sum_{j}\mu(Q_{j})\left(\aver{Q_{j}}(|b|^{\B}+|\nabla b|^{\B})d\mu\right)^{\frac{1}{\B}}
 \\
 &\leq 2CN \|b\|_{W^{1,\B}}\mu(Q_{0})^{1-\frac{1}{\B}}
 \\
 &\leq C\mu(Q_0).
 \end{align*}
 Consequently, the sum in (\ref{bhg}) converges in $W^{1,1}$, $\int g_{0}d\mu=0$ since $\int bd\mu=0$ and $\int h_jd\mu=0$. It follows that $a_0\equiv \frac{g_0}{NCK\alpha \mu(Q_{0})}$ is an $\infty$-atom. Thus we can write
 $$
 b=NCK\alpha \mu(Q_0)a_0+\sum_{j\in \mathbb{N}}h_j.
 $$
 Properties (a) and (d) are then established in this case when $n=1$.
 Property (c) follows from the Whitney covering theorem, since $M$ satisfies $(D)$. We have
 \begin{align*}
 |h_{j}(x)|&\leq |b(x)|+\frac{\mu(Q_j)}{\chi_{j}(Q_j)}\aver{Q_j}|b|d\mu
 \\
 &\leq |b(x)|+ C_0 4^{\frac{D}{q}}K\alpha,
 \end{align*}
 where $D=log_2C_d$ and $C_d$ is the doubling constant. We refer the reader to the proof of Proposition \ref{CZ} for the construction of $h_j$'s and $\chi_j$'s. We have
 $$ h_j=\left(b-\frac{1}{\chi_j(Q_j)} \int_{Q_j} b\chi_j d\mu \right) \chi_j,$$
with $\chi_j(Q_j)=\int_{Q_j} \chi_j d\mu$ and essentially, $\chi_j$ is a smooth version of ${\bf 1}_{Q_j}$, with $|\nabla \chi_j|\leq r_{Q_j}^{-1}$. \\
For $\nabla h_j$ we have
\begin{equation*}
\nabla h_j=\chi_{j} \nabla b +\left(b-\frac{1}{\chi_{j}(Q_j)}\int_{Q_j}b\chi_j d\mu\right)\nabla \chi_j=I+II.
\end{equation*}
We have $|I|\leq |\nabla b|\leq  \mathcal{M}_{q-\epsilon}(\nabla b)$. It remains to estimate II. For $y\in Q_j$
, we have
\begin{align*}
|b(y)-\frac{1}{\chi_{j}(Q_j)}\int_{Q_j}b\chi_j d\mu| & \\
& \hspace{-2cm} \leq \sum_{k=-\infty}^{-1}|b_{Q(y,2^{k}r_{j})}-b_{Q(y,2^{k+1}r_{j})}|d\mu+\left|b_{Q(y,r_{j})}-\frac{1}{\chi_{j}(Q_j)}\int_{Q_j}b\chi_j d\mu\right| \\
& \hspace{-2cm} \leq \sum_{k=-\infty}^{-1}\aver{Q(y,2^{k}r_{j})}|b-b_{Q(y,2^{k+1}r_{j})}| d\mu
+\left|b_{Q(y,r_{j})}-\frac{1}{\chi_{j}(Q_j)}\int_{Q_j}b\chi_jd\mu \right| \\
& \hspace{-2cm} \leq \sum_{k=-\infty}^{-1}\frac{\mu(Q(y,2^{k+1}r_{j}))}{\mu(Q(y,2^{k}r_{j}))}\aver{Q(y,2^{k+1}r_{j})}|b-b_{Q(y,2^{k+1}r_{j})}|d\mu \\
& \hspace{-1cm} +|b_{Q(y,r_j)}-b_{2Q_j}|+\left|\frac{1}{\chi_j(Q_j)}\int_{Q_j}\left( b-\frac{1}{\mu(2Q_j)}\int_{2Q_j}b\right)\chi_j d\mu\right| \\
& \hspace{-2cm} \leq 2^{D}Cr_{j} \mathcal{M}_{q-\epsilon}(|\nabla b|)(y)\sum_{k=-\infty}^{-1}2^{k}+C2^{D/q}3^DK\alpha r_j \\
&  \hspace{-1cm} +\frac{1}{\chi_j(Q_j)}\int_{2Q_j} \left|b-\frac{1}{\mu(2Q_j)}\int_{2Q_j}b d\mu\right||\chi_j| d\mu \\
& \hspace{-2cm} \leq 2^DC r_{j} \mathcal{M}_{q-\epsilon}(|\nabla b|)(y)+Cr_j\left(\aver{Q_j}|\nabla b|^{q-\epsilon}d\mu\right)^{\frac{1}{q-\epsilon}} +C_0 C2^{D/q}K\alpha r_j \\
& \hspace{-2cm} \leq C'\left( \mathcal{M}_{q-\epsilon}(|\nabla b|)(y)+K\alpha\right) r_{j}
\end{align*}
where $C'=\max(C_04^{D/q}, C3^D2^{D/q}, C_0C2^{D/q})$. Thus $|\nabla h_j|\leq  (C'+2) \mathcal{M}_{q-\epsilon}(|\nabla b|)+2C'K\alpha \ind_{Q_j}$. We choose $C''=2C'\geq 1$, and thus (e)
and (f) are proved.
Similarly to (\ref{ebaze}), we deduce (g) and finally property (b) is satisfied by the Whitney covering. The induction hypothesis is then satisfied for $n=1$. \\
We assume that it holds for $n$ and show its validity for $n+1$.
 Consider the set
$$\widetilde{\Omega_{j_n}}:=\left\{ x\in M;  \mathcal{M}(|h_{j_n}|^{q}+|\nabla h_{j_n}|^{q})(x)>(K\alpha)^{q(n+1)}\right\}.$$ Property (g) for $n$ shows that
 $$C_{1}\aver{Q_{j}}(|h_{j_n}|^q+|\nabla h_{j_n}|^q)d\mu\leq C_{1}C^q (K\alpha)^{nq}<(K\alpha)^{(n+1)q}$$ provided $K\alpha>C_{1}^{\frac{1}{q}}C$ and where $C_1$ is the constant in Lemma \ref{coif}. Then Lemma \ref{coif} asserts that $\widetilde{\Omega_{j_n}}\subset 2Q_{j_n}$. Let now $(Q_{j_n,i})$ be a Whitney covering for $\widetilde{\Omega_{j_n}}$.
 We have $\bigcup_{i}Q_{j_n,i}=\widetilde{\Omega_{j_{n}}}\subset \Omega_{n}$ and the $(Q_{j_n,i})_i$ have the bounded overlap property.
  From (c) for $n$, we know that the $(Q_{j_n})$ are $N^{n}$ disjoint. Consequently, the balls $(Q_{j_n,i})$ are $N^{n+1}$ disjoint and therefore we obtain (c) for $n+1$.
 Pose
$$ h_{j_n,i}(x)= \left(h_{j_n}(x)-\frac{1}{\chi_{j_n}^{i}(Q_{j_n,i})}\int_{Q_{j_n,i}}h_{j_n}\chi_{j_n}^{i}d\mu\right)\chi_{j_n}^{i}(x) $$
and $g_{j_n}=h_{j_n}-\sum_{i}h_{j_n,i}$. The same arguments as in Proposition \ref{CZ} show that
$$\|g_{j_n}\|_{W^{1,\infty}}\leq C(K\alpha)^{(n+1)}.$$ Since the support of $h_{j_n}$ is contained in $Q_{j_n}\subset
2Q_{j_n}$ and $\widetilde{\Omega_{j_n}}\subset 2Q_{j_n}$, we deduce that $\supp\,g_{j_n}\subset 2Q_{j_n}$. For every
$i$, $\int h_{j_n,i}\,d\mu=0$  so that (d) follows for $n+1$. We also obtain $\sum_{i}\|h_{j_n,i}\|_{L^1}+\|\,|\nabla
h_{j_n,i} |\,\|_{L^1}\leq C\|h_{j_n}\|_{W^{1,1}}$ as in Proposition \ref{CZ}. Therefore, the equality
$$h_{j_n}=g_{j_n}+\sum_{i}h_{j_n,i}
$$
holds in $W^{1,1}$ and also $\mu-a.e.$ since for each $x$ the sum has at most $N^{(n+1)}$ terms and $\int g_{j_n}d\mu=0$. It
follows that
$$
a_{j_n}=\frac{g_{j_n}}{CN(K\alpha)^{n+1}\mu(2Q_{j_{n}})}$$ is an $\infty$ atom with $\supp \, a_{j_n}\subset 2Q_{j_n}$. We
deduce that the representation (\ref{bex}) holds for $n+1$ and also (a). Let us prove (e) and (f) for $n+1$. The
definition of $h_{j_n,i}$ and (e) for $n$ yield
\begin{align*}
|h_{j_n,i}(x)|&\leq \left(|h_{j_n}(x)|+C
_0\left[\aver{Q_{j_n,i}}|h_{j_n}|^{q}d\mu\right]^{\frac{1}{q}}\right)\chi_{j_n}^{i}(x)
\\
&\leq \left(|b(x)|+2C_04^{\frac{D}{q}}(K\alpha)^{n}+C_04^{\frac{D}{q}}(K\alpha)^{n+1}\right)\chi_{j_n}^{i}(x)
\\
&\leq |b(x)|+2C_04^{\frac{D}{q}}(K\alpha)^{n+1}\chi_{j_n}^{i}(x)
\\
&=|b(x)|+2C'(K\alpha)^{n+1}\chi_{j_n}^{i}(x)
\end{align*}
as long as $K\alpha>2$.
The definition of $\nabla h_{j_n,i}$ and (f) for $n$ yield
\begin{align*}
|\nabla h_{j_n,i}|& \leq |\nabla h_{j_n}|+\left(h_{j_n}-\frac{1}{\chi_{j_n}^i (Q_{j_n,i})}\int_{Q_{j_n,i}}h_{j_n}\chi_{j_n}^i d\mu\right)\nabla \chi_{j_n}^i
\\
& \hspace{-1cm} \leq  (C''+2) \left[K^{n-1}  \mathcal{M}_{q-\epsilon}^n(|\nabla b|)+ (K\alpha)^{n}\ind_{Q_{j_n}}\right]{\bf 1}_{Q_{j_n,i}}   \\
 & \hspace{-0.5cm} +C' \left[ \mathcal{M}_{q-\epsilon}(|\nabla h_{j_n}|) + (K\alpha)^{n+1}{\bf 1}_{Q_{j_n}}\right]{\bf 1}_{Q_{j_n,i}} \\
 & \hspace{-1cm} \leq  (C''+2) \left[K^{n-1}  \mathcal{M}_{q-\epsilon}^n(|\nabla b|)(x)+ (K\alpha)^{n}\ind_{Q_{j_n}}\right]{\bf 1}_{Q_{j_n,i}} \\
  &\hspace{-0.5cm} +C' (C''+2) \left[K^{n-1}  \mathcal{M}_{q-\epsilon}^{n+1}(|\nabla b|)(x) + (K\alpha)^{n}\right]{\bf 1}_{Q_{j_n,i}} +C'(K\alpha)^{n+1}{\bf 1}_{Q_{j_n}}{\bf 1}_{Q_{j_n,i}}\\
 & \hspace{-1cm} \leq (C''+2) \left[K^{n}  \mathcal{M}_{q-\epsilon}^{n+1}(|\nabla b|)(x)+ (K\alpha)^{n+1}{\bf 1}_{Q_{j_n,i}}\right].
\end{align*}
as long as $K,K\alpha$ are large enough (for example we require $K>>4C'$). Now we can prove (b). From (e) and (f), we
deduce that for $x\in \widetilde{\Omega_{j_n}}$,
\begin{align*}
 (K\alpha)^{(n+1)} & <  \mathcal{M}_{q}(|h_{j_n}|+|\nabla h_{j_n}|)(x) \\
 & < (C''+2)\left[ \mathcal{M}_{q}(b)(x)+K^{n-1} {\mathcal M}_q \mathcal{M}_{q-\epsilon}^n(|\nabla b|)(x)+2(K\alpha)^{n} \right].
\end{align*}
provided $K\alpha$ large enough. Thus if we take $K>>4(C''+2)$, we deduce that
$$ \mathcal{M}_{q}\left[|b|+K^{n-1} \mathcal{M}_{q-\epsilon}^{n}(|\nabla
b|)\right](x)>\frac{K^n\alpha^{(n+1)}}{2}$$ and so as $K>1$ we obtain
$$ \mathcal{M}_{q}\left[|b|+ \mathcal{M}_{q-\epsilon}^{n}(|\nabla b|)\right](x)>\frac{K\alpha^{(n+1)}}{2}.$$ Thus
 $\bigcup_{j_n,i}Q_{j_{n},i}\subset \bigcup_{j_n}\widetilde{\Omega_{j_n}}\subset \Omega_{n+1}$.
The last point (g) for $n+1$ is obtained as (\ref{ebaze}) in Proposition \ref{CZ}. The proof is therefore complete.
\findem

\mb We finish this subsection describing the homogeneous version of all these results.

 \begin{df} For $1<\B\leq \infty$, we say that a function $b$ is a homogeneous $(1,\B)$-atom associated to a ball $Q$ if
\begin{itemize}
\item[1.] $b$ is supported in the ball $Q$,
\item[2.] $\|b\|_{\dot{W}^{1,\B}}:=\|\,|\nabla b|\,\|_{L^\B}\leq \mu(Q)^{-\frac{1}{\B'}}$,
\item[3.] $\int b d\mu=0$.
\end{itemize}
\end{df}

\begin{df} For $1<\B\leq \infty$, we define the homogeneous atomic Hardy-Sobolev space $\dot{HS}^{1}_{(\B),ato}$ as follows: $f\in \dot{HS}^{1}_{(\B),ato}$ if $f \in L^{1}_{loc}$ and there exists $(b_{i})_{i}$ a family of homogeneous $(1,\B)$-atoms such that $f=\sum_{i}\lambda_{i}b_{i}$ with $\sum_{i}|\lambda_{i}|<\infty$.
We equip this space with the  semi-norm
$$
\|f\|_{\dot{HS}^{1}_{(\B),ato}}=\inf_{(\lambda_{i})_{i}}\sum_{i}|\lambda_{i}|.
$$
\end{df}
 \begin{prop} $\dot{HS}^{1}_{(\B),ato}/\mathbb{R}$ is a Banach space for every $1<\B\leq \infty$.
\end{prop}
\begin{rem}
Note that every homogeneous $(1,\B)$- atom is an homogeneous $(1,\B')$- atom for $1<\B'\leq \B\leq\infty$ and therefore $\dot{HS}_{(\B),ato}^{1}\subset \dot{HS}_{(\B'),ato}^{1}$ with $\|f\|_{\dot{HS}^{1}_{(\B'),ato}}\leq \|f\|_{\dot{HS}_{(\B),ato}^{1}}$.
\end{rem}
 \begin{prop}\label{DSH} For  $1<\B<\infty$
 the finite subspace $\dot{HS}_{F,(\B),ato}^1$ is dense in $\dot{HS}_{(\B),ato}^1$.
 \end{prop}
  \begin{df} Let $M$ be a Riemmanian manifold. The classical homogeneous Hardy-Sobolev space $\dot{HS}^{1}(M)$ is defined as
$\dot{HS}^{1}=\left\lbrace f\in L^1_{loc}; \nabla f\in H^{1}(M)\right\rbrace$ (see \cite{strichartz2}, \cite{cho} in
the Euclidean case).
 \end{df}
 \begin{prop} $\dot{HS}^{1}$ is a Banach space.
\end{prop}
\begin{prop} We have $\dot{HS}_{(\B),ato}^{1}\subset \dot{HS}^1$ for all $1<\B<\infty$.
\end{prop}

\mb \textbf{Proof of Theorem  \ref{thm:comph}:} Same proof as that of Theorem \ref{thm:comp} but considering the
homogeneous version of the  Calder\'on-Zygmund decomposition.
 \findem



\section{Interpolation of Hardy-Sobolev spaces.} \label{section3}

This section is dedicated to the study of real interpolation of Hardy-Sobolev spaces with Sobolev spaces. First we show how we can use the Calder\'on-Zygmund decomposition for Sobolev functions to obtain interpolation results for the particular Hardy-Sobolev spaces (studied in Subsection \ref{sub:par}). \\
Unfortunately, this method is very specific to this kind of spaces and seems not to be generalized for the study of other Hardy-Sobolev spaces. That is why in Subsection \ref{sub:intabs}, we will use the maximal characterization and the results of Subsection \ref{subsec:maximal} to obtain interpolation results in a  more abstract background.

\subsection{Interpolation of particular Hardy-Sobolev spaces.} \label{sub:intpar}

First as done in \cite{Nadine} and \cite{Nadine1}, we want to prove interpolation results using an adapted ``Calder\'on-Zygmund'' decomposition for Sobolev functions.

\mb Let us describe it~:

\begin{prop}[Calder\'{o}n-Zygmund lemma for Sobolev functions]\label{CZ} Let $M$ be a complete non-compact Riemannian manifold satisfying $(D)$.
Let $1<q<\infty$ and assume that $M$ satisfies a Poincar\'e inequality $(P_{q})$. Let $q\leq p<\infty$, $f \in W^{1,p}$ and $\alpha>0$. Then one can find a collection of balls $(Q_{i})_{i}$, functions $b_{i}$ and a Lipschitz function $g$ such that the following properties hold:
\begin{equation}
f = g+\sum_{i}b_{i} \label{dfaze}
\end{equation}
\begin{equation}
|g(x)|\leq C\alpha\,\textrm{ and }\,|\nabla g(x)|\leq C\alpha\quad \mu-a.e\; x\in M \label{egaze}
\end{equation}
\begin{equation}\label{ebaze}
\supp b_{i}\subset Q_{i}, \|b_{i}\|_{HS_{(q),ato}^{1}}\leq C\alpha\mu(Q_{i})
\end{equation}
\begin{equation}\label{sbaze}
\sum_{i}\mu(Q_{i})\leq C\alpha^{-p}\int (|f|+|\nabla f|)^{p} d\mu
\end{equation}
\begin{equation}\label{rbaze}
\sum_{i}{\bf 1}_{Q_{i}}\leq N
\end{equation}
where $C$ and $N$  only depend on $q$, $p$ and on the constants in $(D)$ and $(P_{q})$.
\end{prop}

\mb This proposition is very similar to the ones of \cite{Nadine,Nadine1}. So we do not detail the proof and just explain the modifications. The new and important fact, is that the functions $b_i$ (appearing in the decomposition) belong to the atomic Hardy-Sobolev spaces and not just to the Sobolev space $W^{1,1}$.

\dem
Let  $f\in W^{1,p}$, $\alpha>0$ and consider
$\Omega=\left\lbrace x \in M:  \mathcal{M}(|f|^q+|\nabla f|^q)(x)>\alpha^{q}\right\rbrace$. If $\Omega=\emptyset$, then set
$$
 g=f\;,\quad b_{i}=0 \, \text{ for all } i
$$
so that (\ref{egaze}) is satisfied according to the Lebesgue differentiation theorem. Otherwise, the maximal theorem yields
\begin{align}
    \mu(\Omega)&\leq C\alpha^{-p}\|(|f|+ |\nabla f|)^{q}\|_{\frac{p}{q}}^{\frac{p}{q}} \nonumber\\
            & \leq C \alpha^{-p} \Bigr(\int | f|^{p} d\mu +\int |\nabla f|^{p} d\mu\Bigl) \label{mOaze}
\\
            &<+\infty. \nonumber
\end{align}
 In particular $\Omega \neq M$ as $\mu(M)=+\infty$. Let $F$ be the complement of $\Omega$. Since $\Omega$ is an open set distinct of $M$, let
$(\underline{Q_{i}})$ be a Whitney decomposition of $\Omega$ (\cite{coifman1}).
That is,  the $\underline{Q_{i}}$ are pairwise disjoint, and there exist two constants $C_{2}>C_{1}>1$, depending only on the metric, such that
\begin{itemize}
\item[1.] $\Omega=\cup_{i}Q_{i}$ with $Q_{i}=
C_{1}\underline{Q_{i}}$ and the balls $Q_{i}$ have the bounded overlap property;
\item[2.] $r_{i}=r(Q_{i})=\frac{1}{2}d(x_{i},F)$ and $x_{i}$ is
the center of $Q_{i}$;
\item[3.] each ball $\overline{Q_{i}}=C_{2}Q_{i}$ intersects $F$ ($C_{2}=4C_{1}$ works).
\end{itemize}
For $x\in \Omega$, denote $I_{x}=\left\lbrace i:x\in Q_{i}\right\rbrace$.
Recall that $\sharp I_{x} \leq N$ and fixing $j\in I_{x}$, $Q_{i}\subset 7Q_{j}$ for all $i\in I_{x}$. \\
Conditions (\ref{rbaze}) and (\ref{sbaze}) are satisfied due to (\ref{mOaze}). Using the doubling property, we have
\begin{equation}\label{faze}
\int_{Q_{i}} (|f|^{q}+|\nabla f|^{q})d\mu \leq C \alpha^{q}\mu(Q_{i}).
\end{equation}
Let us now define the functions $b_{i}$. For this, we construct a partition of unity $(\chi_{i})_{i}$ be  a partition of unity of $\Omega$ subordinated to the covering $(Q_{i})$. Each $\chi_{i}$ is a Lipschitz function supported in $Q_{i}$ with
$\displaystyle\|\,|\nabla \chi_{i}|\, \|_{\infty}\leq
\frac{C}{r_{i}}$.\\
We set $b_{i}=(f-\frac{1}{\chi_{i}(Q_{i})}\int_{Q_{i}}f\chi_{i})\chi_{i}$ where $\chi_{i}(Q_{i})\simeq \mu(Q_i)$ means $\int_{Q_{i}}\chi_{i}d\mu$. This is the main change, which is necessary as we look for a vanishing mean value for $b_i$'s. \\
By usual arguments and Poincar\'e inequality $(P_q)$, we can estimate $b_i$ in the Sobolev space $W^{1,q}$: $\|b_i\|_{W^{1,q}}\leq C\alpha \mu(Q_i)^{\frac{1}{q}}$.
Then by writing $b_{i}=\alpha_{i}\alpha_{i}^{-1}b_{i}=\alpha_{i}a_{i}$ with $\alpha_{i}=C\alpha\mu(Q_{i})$, we deduce that the functions $a_{i}$ are $(1,q)$-atoms -- and in fact $(1,r)$-atoms for every $r\leq q$-- associated to the ball $Q_i$. Therefore $b_{i}\in HS_{(q),ato}^{1}$ with $\|b_{i}\|_{HS_{(q),ato}^{1}}\leq \alpha_{i}=C\alpha \mu(Q_{i})$ and also $b=\sum_{i}b_{i}\in HS_{(q),ato}^{1}$ with $\|b\|_{HS_{(q),ato}^{1}}\lesssim \alpha \mu(\Omega)$.
Thus (\ref{ebaze}) is proved.
Set $\displaystyle g=f-\sum_{i}b_{i}$. Since the sum is locally finite on $\Omega$, as usually $g$ is defined  almost everywhere on $M$ and $g=f$ on $F$. Moreover, $g$ is a locally integrable function on $M$.
It remains to prove (\ref{egaze}). We have
\begin{align*}
\nabla g &= \nabla f -\sum_{i}\nabla b_{i}
\\
&=\nabla f-(\sum_{i}\chi_{i})\nabla f -\sum_{i}(f-\frac{1}{\chi_{i}(Q_{i})}\int_{Q_{i}}f\chi_{i}d\mu)\nabla
\chi_{i}
\\
&=\ind_{F}(\nabla f) - \sum_{i}(f-\frac{1}{\chi_{i}(Q_{i})}\int_{Q_{i}}f\chi_{i}d\mu)\nabla
\chi_{i}.
\end{align*}
From the definition of $F$ and the Lebesgue differentiation theorem, we have $\ind_{F}(|f|+|\nabla f|)\leq \alpha\;\mu -$a.e. We claim that a similar estimate holds for $$h=\sum_{i}\left(f-\frac{1}{\chi_{i}(Q_{i})}\int_{Q_{i}}f\chi_{i}d\mu\right)\nabla
\chi_{i},$$
 that is $|h(x)|\leq C\alpha$ for all $x\in M$. For this, note first that $h$ vanishes on $F$ and the sum defining $h$ is locally finite on $\Omega$.
Then fix  $x\in \Omega$ and $j \in I_x$. Note that $\displaystyle \sum_{i}\chi_{i}(x)=1$ and $\displaystyle \sum_{i}\nabla\chi_{i}(x)=0$, so
$$h(x)=\sum_{i\in I_x} \left[\left(\frac{1}{\mu(7Q_{j})}\int_{7Q_{j}}f d\mu\right) -\left(\frac{1}{\chi_{i}(Q_{i})}\int_{Q_{i}}f\chi_{i}d\mu\right)\right]\nabla
\chi_{i}(x).$$
For all $i,j\in I_{x}$, by the construction of  the Whitney collection, the ball $Q_i$ and $Q_j$ have equivalent radius and  $Q_i \subset 7Q_j$.
Thus
\begin{align}
\left|\frac{1}{\chi_{i}(Q_{i})}\int_{Q_{i}}f\chi_{i}d\mu-\aver{7Q_{j}}f
d\mu \right| & \nonumber \\
 & \hspace{-3cm} \leq \frac{1}{\chi_{i}(Q_{i})} \int_{Q_{i}} \left|f - \aver{7Q_{j}}fd\mu \right| |\chi_i| d\mu  \nonumber \\
 & \hspace{-3cm} \lesssim  \aver{7Q_{j}} \left|f-\aver{7Q_{j}}fd\mu \right|  d\mu \nonumber \\
 & \hspace{-3cm} \lesssim r_j \left( \aver{7Q_j} |\nabla f|^q d\mu \right)^{1/q} \nonumber \\
 & \hspace{-3cm} \lesssim \alpha r_j. \label{relim}
 \end{align}
We used $(D)$, $(P_{q})$, $\chi_i(Q_i)\simeq \mu(Q_i)$ and (\ref{faze})  for $7Q_j$. Hence
\begin{align}
|h(x)| \lesssim \sum_{i\in I_{x}}\alpha r_j r_{j}^{-1} \leq CN\alpha .
\end{align}
Then the end of the proof is classical and is exactly the same as that of  the decompositions proved in \cite{Nadine,Nadine1}. We do not repeat it. \findem

\mb According to \cite{Nadine,Nadine1}, we know how to obtain interpolation results from an adapted ``Calder\'on-Zygmund decomposition''. We quickly recall them (for an easy reference) in order to obtain a real interpolation result between the Hardy-Sobolev spaces $HS_{(q),ato}^1$ and Sobolev spaces.

\mb First we characterize the $K$-functional of real interpolation in the following theorem:

\begin{prop} \label{EKHS} Let $M$ be a complete Riemannian manifold satisfying the doubling condition $(D)$ and Poincar\'e inequality $(P_q)$ for some $q\in(1,\infty)$. Then
\begin{itemize}
\item[1.] for all  $r\in(1,\infty)$, there exists $C_{1}>0$ such that for every $f\in HS_{(r),ato}^{1}+W^{1,\infty}$ and $t>0$,
$$
K(f,t, HS_{(r),ato}^{1},W^{1,\infty})\geq C_{1}t\left(|f|^{**}+|\nabla f|^{**}\right)(t);
$$
\item[2.] for $1<q\leq p<\infty$, there exists $C_{2}>0$ such that for every $f\in W^{1,p}$ and $t>0$,
$$
K(f,t, HS_{(r),ato}^{1},W^{1,\infty})\leq C_{2}t\left(|f|^{q**\frac{1}{q}}+|\nabla f|^{q**\frac{1}{q}}\right)(t).
$$
\end{itemize}
We have the same results replacing the space $HS_{(r),ato}^{1}$ by $HS^{1}$.
\end{prop}

\dem We only write the proof for the space $HS_{(r),ato}^{1}$. We have already seen (Theorem \ref{thm:comp}) that under our assumption $HS_{(r),ato}^{1}=HS_{(q),ato}^{1}$ for $r\in [q,\infty]$. We just have to prove our result for $r\in(1,q]$. The lower bound of $K$ is trivial. It follows from the characterization of $K$ between  $L^1$ and $L^\infty$.  Now for the upper bound of $K$ of point 2., take  $f\in W^{1,p}$ and $q\leq p<\infty$.
Let $t>0$. We consider the Calder\'{o}n-Zygmund decomposition of Proposition \ref{CZ} for $f$ with
$\alpha=\alpha(t)=\left( \mathcal{M}(|f|^q+|\nabla f|^q)\right)^{*\frac{1}{q}}(t)$.
We write $ \displaystyle f=\sum_{i}b_{i}+g=b+g $ where
$(b_{i})_{i},\,g$ satisfy the properties of the proposition. From the bounded overlap property of the $B_{i}$'s, it follows that for all $r\leq q$
\begin{align*}
\| b \|_{HS_{(r),ato}^{1}}&\leq N \sum_{i}
\|b_{i}\|_{HS_{(r),ato}^{1}}
\\
&\leq C\alpha^{r}(t)\sum_{i}\mu(Q_{i})
\\
&\leq C\alpha^{r}(t)\mu(\Omega).
\end{align*}
Moreover, since $( \mathcal{M}f)^{*}\sim f^{**}$ and $(f+g)^{**}\leq f^{**}+g^{**}$ (c.f \cite{bennett},\cite{bergh}) , we get $$
\alpha(t)\lesssim \left(|f|^{q**{\frac{1}{q}}}(t)+|\nabla
f|^{q**{\frac{1}{q}}}(t)\right).
$$
Noting that for this choice of $\alpha(t)$, $\mu(\Omega_t)\leq t$ (c.f \cite{bennett},\cite{bergh}), we deduce that
 \begin{align}
 K(f,t,HS_{(r),ato}^{1},W^{1,\infty}) & \leq \|b\|_{HS^{1}_{(r),ato}} + t\|g\|_{W^{1,\infty}} \nonumber \\
 & \lesssim t\left(|f|^{q**\frac{1}{q}}(t)+|\nabla
f|^{q**\frac{1}{q}}(t)\right) \label{Kr}
\end{align}
for all $t>0$ and obtain the desired inequality for $f\in W^{1,p},\, q\leq p<\infty$. \findem

\mb Then integrating the $K$-functional yields

\begin{prop} \label{IHS} Let $M$ be a complete Riemannian manifold satisfying $(D)$ and $(P_q)$, for some $1< q <\infty$. Then for all $r\in(1,\infty]$ and $p\in(q,\infty)$, $W^{1,p}$ is a real interpolation space between $HS_{(r),ato}^{1}$ and $W^{1,\infty}$. More precisely, we have
$$ \left(HS_{(r),ato}^{1},W^{1,\infty}\right)_{1-\frac{1}{p},p} = W^{1,p}.$$
\end{prop}

\mb We refer the reader to the previously cited papers for a detailled proof. We also have an analogous interpolation result for the Hardy-Sobolev space $HS^{1}$ instead of $HS_{(r),ato}^{1}$. Note that $HS^{1}_{(r),ato}\subset HS^1$ and $\|f\|_{HS^1}\leq 2\|f\|_{HS^1_{(r),ato}}$.
 \\
 \\
\textbf{Proof of Theorem \ref{thm:hardy}:} The proof follows from Proposition \ref{IHS} and the Reiteration Theorem
(see \cite{bennett}, Theorem 2.4).
 \findem

\mb All these results are based on the well adapted ``Calder\'on-Zygmund decomposition''. The first one (described in
\cite{A1} by P. Auscher) was written for homogeneous Sobolev spaces. We can write an analog result of Proposition
\ref{CZ} for homogeneous Sobolev spaces. Then we estimate the functional $K$ (as in \cite{Nadine}) and obtain the
homogeneous interpolation Theorem \ref{thm:hardyhomo}:

\gb \textbf{Proof of Theorem \ref{thm:hardyhomo}: } Analogous proof to that of Theorem  \ref{thm:hardy} and \ref{IHS}.
\findem

\mb We used a ``Calder\'on-Zygmund'' decomposition to obtain an interpolation result for the particular Hardy-Sobolev spaces. These arguments give positive interpolation results under the assumptions of doubling property and Poincar\'e inequality. Unfortunately, this method seems not to work for abstract and more general Hardy-Sobolev spaces: the way to make appear the ``atoms'' is very particular. That is why, in the next subsection, we develop other arguments to obtain interpolation results with abstract Hardy-Sobolev spaces. We will use our maximal characterization of Sobolev spaces (Subsection \ref{subsec:maximal}) and ideas of \cite{BJ}.

\subsection{Interpolation of abstract Hardy-Sobolev spaces.} \label{sub:intabs}

\mb We refer the reader to Subsection \ref{subsec:method1} for the definition of abstract Hardy-Sobolev spaces associated to a collection of ``local operators'' ${\mathbb B}$.

\mb To prove our results, we will follow ideas of \cite{BJ} and \cite{B2} using duality and some maximal operators associated to the collection ${\mathbb B}$. Let us first define them.

\begin{df} Let $\sigma\in(1,\infty]$. We set $A_Q=Id-B_Q$ and
\be{opeM} \forall x\in M, \qquad M_{{\mathbb B},\sigma}(f)(x):= \sup_{\genfrac{}{}{0pt}{}{Q \textrm{ball}}{x\in Q}}\  \frac{1}{\mu(Q)^{1/\sigma}}  \left\|A_Q^*(f)\right\|_{W^{-1,\sigma}(Q)}. \ee
We define a sharp maximal function adapted to our operators. For $s>0$,
$$ \forall x\in M, \qquad M^\sharp_{{\mathbb B},s} (f)(x):= \sup_{\genfrac{}{}{0pt}{}{Q \textrm{ball}}{x\in Q}} \frac{1}{\mu(Q)^{1/s}} \left\|B_Q^* (f)\right\|_{W^{-1,s}(Q)}.$$
\end{df}

\mb We refer the reader to Definition \ref{ensembleI} for the notation ${\mathcal I}_M$ and Subsection \ref{subsec:maximal} for the definition of some maximal operators and the assumption (\ref{Hp}) . \\
We can now prove Theorem  \ref{thm:generalbis} .
\begin{rem}
We want to emphasize that we only require the use of the ``finite Hardy-Sobolev'' space $HW^{1}_{F,ato}$. With our new
maximal operators, the assumption (\ref{assumpintro}) can be written as \be{assu} M_{{\mathbb B},\sigma} \lesssim
M_{S,*,\B'}.\ee
\end{rem}
\textbf{Proof of Theorem \ref{thm:generalbis}:}
 From the $HW^1_{F,ato}-L^1$ boundedness, it is quite easy to check that for each ball $Q$, the operator $TB_Q$ is bounded from $W^{1,\B}(Q)$ into $L^1$ with
$$ \| TB_Q\|_{W^{1,\B}(Q)\to L^1} \lesssim \mu(Q)^{1/\B'}.$$
By duality, we deduce that $B_Q^*T^*$ is bounded from $L^\infty$ to $W^{-1,\B'}$ with
$$ \|B_Q^* T^*\|_{L^\infty\to W^{-1,\B'}(Q)} \lesssim \mu(Q)^{1/\B'}.$$
Thus, we obtain the first inequality~ \be{rel1} \forall f\in L^\infty, \qquad  \left\| M^\sharp_{{\mathbb B},\B'}(T^*f)
\right\|_{L^\infty} \lesssim \|f\|_{L^\infty}. \ee
Now using (\ref{assu}), we obtain
\begin{equation}\label{sharp}
M^\sharp_{{\mathbb B},\B'} \leq M_{S,\B'}+M_{{\mathbb B},\B'} \leq M_{S,\B'}+M_{{\mathbb B},\sigma} \lesssim M_{S,*,\B'}\lesssim M_{S,*,p_0'}.
\end{equation}
Then (\ref{sharp}) with Proposition \ref{prop:importante}, yields the following ``weak type inequality''~
\be{rel2} \forall f\in L^{p_0'}, \qquad  \left\| M^\sharp_{{\mathbb B},\B'}(T^*f) \right\|_{L^{p_0',\infty}} \lesssim \|T^*f\|_{W^{-1,p_0'}} \lesssim \|f\|_{L^{p_0'}}. \ee
Interpolating (\ref{rel1}) and (\ref{rel2}) gives
\be{rel3} \forall q\in(p_0',\infty),\;\forall f\in L^{p_0'}\cap L^q, \qquad \left\| M^\sharp_{{\mathbb B},\B'}(T^*f) \right\|_{L^q} \leq c \|f\|_{L^q}. \ee
Now we use a ``good lambdas'' argument to compare the maximal operators. We use a Sobolev-version of the result of P. Auscher and J.M. Martell: \cite{AM}, Theorem 3.1. With its notation, take a function $F$. We define for all balls $Q$
$$ G_Q= B_Q^*F \textrm{    and    } H_Q= A_Q^*F.$$
The assumption (\ref{assu}) shows that \be{hyp1} \mu(Q)^{-1/\sigma} \left\| H_Q \right\|_{W^{-1,\sigma}} \lesssim
M_{{\mathbb B},\sigma}(F)\lesssim M_{S,*,\B'}(F).\ee By definition of $ M^\sharp_{{\mathbb B},\B'} $, we have \be{hyp2}
\mu(Q)^{-1/\B'}\left\|G_Q \right \|_{W^{-1,\B'}} \lesssim M^\sharp_{S,\B'}(F). \ee From these two inequalities, we
claim that  the following good lambda inequality holds (for $K$ large enough and $\gamma$ as small as we want)~
\be{googl} \mu\left( \left\{ M_{S,*,\B'}(F)>K\lambda,\ M^\sharp_{{\mathbb B},\B'}(F)\leq \gamma\lambda \right\}\right)
\lesssim (K^{-\sigma}+\gamma^{\B'} K^{-\B'}) \mu\left(\left\{M_{S,*,\B'}(F) >\lambda \right\}\right). \ee
We postpone the proof of this claim to Lemma \ref{lemmag}. As usually this inequality is satisfied for all $\lambda>0$ if $\mu(X)=\infty$ and only for $\lambda \gtrsim \|M_{S,*,\B'}(F)\|_{L^1}$ if the measure is finite. \\
Assuming this fact, we will conclude the proof. By classical arguments (see proof of Theorem 3.1 in \cite{AM}) we deduce that for $p_0\in(\sigma',\B)$ if $M_{S,*,\B'}(F)\in L^{p_0',\infty}$ then
$$ \| M_{S,*,\B'}(F) \|_{L^q} \lesssim \| M^\sharp_{{\mathbb B},\B'}(F) \|_{L^q} + \| M_{S,*,\B'}(F) \|_{L^1} {\bf 1}_{\mu(X)<\infty}$$
for all $q\in(p_0',\sigma)$ with an implicit constant depending on $q$.
Now we take a function $h \in L^{p_0'} \cap L^q$. Denoting $F=T^*(h)$, we have $F\in W^{-1,p_0'}$. Proposition  \ref{prop:importante} shows that $M_{S,*,\B'}(F)$ belongs to $L^{p_0',\infty}$. Thus we can apply the previous inequality which together with (\ref{rel3}) yield
\begin{align*}
 \| M_{S,*,\B'}(T^* h) \|_{L^q} & \lesssim \| M^\sharp_{{\mathbb B},\B'}(T^* h) \|_{L^q}+\| M_{S,*,\B'}(T^* h) \|_{L^1} {\bf 1}_{\mu(X)<\infty} \\
 & \lesssim \|h\|_{L^q} + \| M_{S,*,\B'}(T^* h) \|_{L^1} {\bf 1}_{\mu(X)<\infty}.
\end{align*}
If the space $X$ is of finite measure, using the $W^{1,p_0}-L^{p_0}$ boundedness of $T$ and Proposition  \ref{prop:importante}, we remark that
$$ \| M_{S,*,\B'}(T^* h) \|_{L^1} \lesssim \| M_{S,*,\B'}(T^* h) \|_{L^{p_0'}} \lesssim \|T^*(h) \|_{W^{-1,p_0'}} \lesssim \|h\|_{L^{p_0'}} \lesssim \|h\|_{L^q}.$$
This inequality with the fact $(\B',q)\in {\mathcal I}_M$ -- since $q'\in (\sigma',p_0)$ --, shows that
$$ \forall h\in L^q \cap L^{p_0}, \qquad \|T^*h\|_{W^{-1,q}} \lesssim \|h\|_{L^q}.$$
By duality, we deduce that there is a constant $c=c(p)$ such that
\begin{equation}\label{T} \forall f\in W^{1,p_0} \cap W^{1,q'}, \qquad \|T(f)\|_{L^q} \leq c \|f\|_{W^{1,q'}}.
\end{equation}
Consequently, inequality (\ref{T}) holds for all $q\in (p'_0,\sigma)$, and therefore $T$ admits a continuous extension from $W^{1,p}$ to $L^p$ for all $p\in
 (\sigma',p_0)$.
\findem

\mb It remains  to prove (\ref{googl}).

\begin{lem} \label{lemmag} With the notations of the previous proof, we have the following good lambda inequality. For all $\lambda>0$ (or only for $\lambda \gtrsim \|M_{S,*,\B'}(F)\|_{L^1}$ if the measure is finite)
$$ \mu\left( \left\{ M_{S,*,\B'}(F)>K\lambda,\ M^\sharp_{{\mathbb B},\B'}(F)\leq \gamma\lambda \right\}\right) \lesssim (K^{-\sigma}+\gamma^{\B'} K^{-\B'}) \mu\left(\left\{M_{S,*,\B'}(F) >\lambda \right\}\right).$$
\end{lem}

\dem The proof is exactly the same as that of Theorem 3.1 in \cite{AM}, adapted to our maximal operators. We deal only with the case when $\mu(X)=\infty$. We consider the sets
$$ B_\lambda:= \left\{ M_{S,*,\B'}(F)>K\lambda,\ M^\sharp_{{\mathbb B},\B'}(F)\leq \gamma\lambda \right\}$$
and
$$ E_\lambda:= \left\{ M_{S,*,\B'}(F)>\lambda \right\}.$$
First since $K\geq 1$, we have $B_\lambda\subset E_\lambda$. We choose $(Q_j)_j$ a Whitney decomposition of $E_\lambda$ and write $x_j$ for a point in $4 Q_j \cap E_\lambda^c$. Let $j$ such that $B_{\lambda}\cap Q_j\neq\emptyset$ and $x\in B_\lambda \cap Q_j$. We have
\be{eq10} M_{S,*,\B'}(F)(x):=\inf_{F=\phi_0-div(\phi_1)} \sup_{\genfrac{}{}{0pt}{}{Q \textrm{ball}}{x\in Q}}\  \frac{1}{\mu(Q)^{1/\B'}}  \left\||\phi_0|+|\phi_1|\right\|_{L^{\B'}(Q)} \geq K\lambda. \ee
Let $F=\psi_0-div(\psi_1)$ and $Q_{ext}$ be an extremize decomposition and ball of (\ref{eq10}).
Assume first that $Q_{ext}$ satisfies $Q_{ext}\cap (8Q_j)^c \neq \emptyset$. Since $x_j\in 4Q_{ext}$ and
$$ \inf_{F=\phi_0-div(\phi_1)} \sup_{Q\ni x_j} \frac{1}{\mu(Q)^{1/\B'}}  \left\||\phi_0|+|\phi_1|\right\|_{L^{\B'}(Q)}\leq M_{S,*,\B'}(F)(x_j)\leq \lambda,$$
we  deduce that
$$ M_{S,*,\B'}(F)(x) \leq \left(\frac{\mu(4Q_{ext})}{\mu(Q_{ext})} \right)^{1/\B'} \lambda.$$
Therefore, for a large enough constant $K$, the doubling property of the measure shows that the assumption $Q_{ext}\cap (8Q_j)^c \neq \emptyset$ is false. We deduce that $Q_{ext} \subset 8Q_j$ and therefore
$$ M_{S,*,\B'}(F)(x)=\inf_{F=\phi_0-div(\phi_1)} \sup_{\genfrac{}{}{0pt}{}{Q \textrm{ball}}{x\in Q\subset 8 Q_j}}\  \frac{1}{\mu(Q)^{1/\B'}}  \left\||\phi_0|+|\phi_1|\right\|_{L^{\B'}(Q)} \geq K\lambda.$$
Write $F= B_{8Q_j}^*F + A_{8Q_j}^*F$. It follows that
\begin{align*}
 \mu(B_\lambda \cap Q_j) \leq & \mu\left(\left\{ \inf_{B_{8Q_j}^*F=\phi_0-div(\phi_1)} \sup_{\genfrac{}{}{0pt}{}{Q \textrm{ball}}{x\in Q\subset 8 Q_j}}\  \frac{1}{\mu(Q)^{1/\B'}}  \left\||\phi_0|+|\phi_1|\right\|_{L^{\B'}(Q)} \geq K\lambda/2\right\} \right) \\
 & \hspace{-1cm} + \mu\left(\left\{ \inf_{A_{8Q_j}^*F=\phi_0-div(\phi_1)} \sup_{\genfrac{}{}{0pt}{}{Q \textrm{ball}}{x\in Q\subset 8 Q_j}}\  \frac{1}{\mu(Q)^{1/\B'}}  \left\||\phi_0|+|\phi_1|\right\|_{L^{\B'}(Q)} \geq K\lambda/2\right\} \right).
\end{align*}
The first term is controlled by the ``weak type $(\B',\B')$'' of the maximal operator $M_{S,*,\B'}$ (local version of Proposition \ref{prop:importante})~:
\begin{align}\label{ee}
 & \mu\left(\left\{ \inf_{B_{8Q_j}^*F=\phi_0-div(\phi_1)} \sup_{\genfrac{}{}{0pt}{}{Q \textrm{ball}}{x\in Q\subset 8 Q_j}}\  \frac{1}{\mu(Q)^{1/\B'}}  \left\||\phi_0|+|\phi_1|\right\|_{L^{\B'}(Q)} \geq K\lambda/2\right\} \right) \hspace{1cm} & \nonumber\\
 & \hspace{6cm} \lesssim \frac{1}{K^{\B'}\lambda^{\B'}} \| M_{S,*,\B'}(B_{8Q_j}^*F)\|_{L^{\B',\infty},8Q_j}^{\B'} &\nonumber \\
 & \hspace{6cm} \lesssim \frac{1}{K^{\B'}\lambda^{\B'}} \| B_{8Q_j}^*F \|_{W^{-1,\B'}(8Q_j)}^{\B'} & \nonumber\\
 & \hspace{6cm} \lesssim \frac{1}{K^{\B'}\lambda^{\B'}} \mu(Q_j) \inf_{8 Q_j} M^\sharp_{{\mathbb B},\B'}(F)^{\B'}  &\nonumber\\
 &\hspace{6cm} \lesssim \frac{\gamma^{\B'}}{K^{\B'}} \mu(Q_j).&
\end{align}
For the last inequality, we used the fact that $B_{\lambda}\cap Q_j\neq \emptyset$.
For the second term, we use similar arguments with $\B'\leq \sigma$
\begin{align*}
& \mu\left(\left\{ \inf_{A_{8Q_j}^*F=\phi_0-div(\phi_1)} \sup_{\genfrac{}{}{0pt}{}{Q \textrm{ball}}{x\in Q\subset 8 Q_j}}\  \frac{1}{\mu(Q)^{1/\B'}}  \left\||\phi_0|+|\phi_1|\right\|_{L^{\B'}(Q)} \geq K\lambda/2\right\} \right) \hspace{2cm} & \\
 & \hspace{8cm} \lesssim \frac{1}{K^\sigma\lambda^\sigma} \|M_{S,*,\sigma} A_{8Q_j}^*F\|_{L^{\sigma,\infty},8Q_j}^\sigma &
 \\
 &\hspace{8cm} \lesssim\frac{1}{K^\sigma\lambda^\sigma} \| A_{8Q_j}^*F\|_{W^{-1,\sigma}(8Q_j)}^\sigma. &
\end{align*}
The above assumption (\ref{assu}) shows that
\begin{align*}
 \| A_{8Q_j}^*F\|_{W^{-1,\sigma}(8Q_j)}^\sigma & \lesssim \mu(Q_j) \inf_{8Q_j} M_{{\mathbb B},\sigma}(F)^\sigma \lesssim \mu(Q_j) \inf_{8Q_j} M_{S,*,\B'}(F)^\sigma \\
 & \lesssim \lambda^\sigma \mu(Q_j).
\end{align*}
We used in the last inequality that $x_j\in 8Q_j$ and $M_{S,*,\B'}(F)(x_j)\leq \lambda$. Thus, we proved an analogous inequality of (\ref{ee}) for the second term. We deduce that
$$ \mu\left(B_\lambda \cap Q_j \right) \lesssim \left(\frac{\gamma^{\B'}}{K^{\B'}} + \frac{1}{K^\sigma}\right) \mu(Q_j). $$
Summing over $j$, the proof is therefore complete. \findem

\mb In the next proposition, we give a useful criterion to insure
 the main assumption (\ref{assu})~:

\begin{prop} \label{prop:utile} Assume that the operators $A_Q$ satisfy
$$\forall j\geq 0 \quad \frac{1}{\mu(2^{j+1}Q)^{1/\B}} \left\| A_Q(f)\right\|_{W^{1,\B}(S_j(Q))} \leq \alpha_j(Q) \frac{1}{\mu(Q)^{1/\sigma'}}\left\|f \right\|_{W^{1,\sigma'}(Q)},$$
for all functions $f$ supported in the ball $Q$, where the coefficients $\alpha_j(Q)$ satisfy
 \be{coef} \sup_{Q \textrm{ ball}}\ \sum_{j\geq 0} \frac{\mu(2^{j+1}Q)}{\mu(Q)} \alpha_j(Q) <\infty. \ee
Then the maximal operator $M_{{\mathbb B},\sigma}$ is bounded by $M_{S,*,\B'}$.
\end{prop}

\dem  Let $x\in M$. For a ball $Q$, we denote $S_j(Q)=2^jQ \setminus 2^{j-1}Q$. We estimate the Sobolev-norm by
duality~
\begin{align*}
M_{{\mathbb B},\sigma} (f)(x) & =\sup_{Q;\,x\in Q} \sup_{\genfrac{}{}{0pt}{}{g \in C^\infty_0(Q)}{\|g\|_{W^{1,\sigma'}}\leq 1}} \mu(Q)^{-1/\sigma} \int A_Q^*(f) g d\mu \\
 & = \sup_{Q;\,x \in Q} \sup_{\genfrac{}{}{0pt}{}{g \in C^\infty_0(Q)}{\|g\|_{W^{1,\sigma'}} \leq 1}} \mu(Q)^{-1/\sigma} \int f A_Q(g) d\mu.
\end{align*}
Take a decomposition $f=\phi_0-div(\psi_0)$. Then we have
\begin{align*}
M_{{\mathbb B},\sigma}(f)(x) & \leq \sup_{Q;\,x \in Q} \sup_{\genfrac{}{}{0pt}{}{g \in C^\infty_0(Q)}{\|g\|_{W^{1,\sigma'}} \leq 1}} \mu(Q)^{-1/\sigma} \sum_{j\geq 0} \int_{S_j(Q)} \left[\phi_0 A_Q(g)+\psi_0 \nabla A_Q(g)\right]  d\mu  \\
 & \leq \sup_{Q;\,x \in Q} \mu(Q)^{-1/\sigma} \sup_{\genfrac{}{}{0pt}{}{g \in C^\infty_0(Q)}{\|g\|_{W^{1,\sigma'}} \leq 1}} \sum_{j\geq 0}  \left\||\phi_0|+|\psi_0| \right\|_{L^{\B'}(S_j(Q))} \left\| A_Q(g) \right\|_{W^{1,\B}(S_j(Q))}.
\end{align*}
Our assumption yields
\begin{align*}
M_{{\mathbb B},\sigma}(f)(x) & \leq \sup_{Q;\,x \in Q} \mu(Q)^{-1/\sigma}  \sum_{j\geq 0}  \left\||\phi_0|+|\psi_0| \right\|_{L^{\B'}(S_j(Q))} \alpha_j(Q) \frac{\mu\left(2^{j+1}Q\right)^{1/\B}}{\mu(Q)^{1/\sigma'}} \\
 & \leq \sup_{Q;\,x \in Q} \sum_{j\geq 0}  \left\||\phi_0|+|\psi_0| \right\|_{L^{\B'}(2^jQ)} \mu(2^{j+1} Q)^{-1/\B'} \alpha_j(Q) \frac{\mu(2^{j+1}Q)}{\mu(Q)} \\
 & \leq M_{HL,\B'}(|\phi_0|+|\psi_0|)(x) \sup_{Q;\,x \in Q}  \sum_{j\geq 0}   \alpha_j(Q) \frac{\mu(2^{j+1}Q)}{\mu(Q)} \\
 & \lesssim M_{HL,\B'}(|\phi_0|+|\psi_0|)(x).
\end{align*}
These inequalities hold for every decomposition $f=\phi-div(\psi)$. Taking the infimum over all these decompositions, we obtain the desired inequality.
\findem

\mb With an extra assumption (as in \cite{B2}), we  obtain the real interpolation result of Theorem \ref{thm:general2} :
 \\
\textbf{Proof of Theorem \ref{thm:general2}:} The proof is the same as the one of Theorem 3.14 in \cite{B2} using the arguments of Theorem \ref{thm:generalbis}. We omit it.
\findem

\mb Let us  compare our assumption $(\B',p_\theta')\in {\mathcal I}_M$ with Poincar\'e inequality~:

\begin{rem}  Assume that $\B'\leq p_\theta'$ (else $(\B',p_\theta')\in {\mathcal I}_M$ is always satisfied, see \cite{BB}) and $p_\theta\leq 2$. Thanks  to Theorem \ref{poincare}, we can check that the assumption $(\B',p_\theta')\in {\mathcal I}_M$ is implied by the Poincar\'e inequality $(P_{p_\theta})$ if $\B'\geq 2$, which corresponds to a variant of the assumption done in \cite{Nadine} (in \cite{Nadine}, the author used local hypotheses of doubling and Poincar\'e, here we are under the global hypotheses) to interpolate the corresponding non-homogeneous Sobolev spaces.
\end{rem}

\section{Applications}\label{section:App}


\subsection{Operators with regularity assumptions about the kernel.}

In this subsection, we look for a ``Sobolev'' version of results for Calder\'on-Zygmund operators on Lebesgue spaces.

\begin{df} Let $T$ be a linear operator bounded from $\dot{W}^{1,p_0}$ (resp. $W^{1,p_0}$) to $L^{p_0}$. We say that it is associated to a kernel $K(x,y)$ if for every compactly supported function $f$ and $x\in supp(f)^c$ we have the integral representation~:
$$ T(f)(x) = \int K(x,y) f(y) d\mu(y).$$
\end{df}

\mb We introduce the following regularity property for such kernel~:

\be{kernel} \Gamma:=\sup_{Q \textrm{ ball}} \sup_{y,z\in Q} r_Q \int_{M\setminus 4Q} \left|K(x,y)-K(x,z) \right| d\mu(x)<\infty. \ee

\mb This subsection is devoted to the study of operators $T$ associated to a kernel satisfying (\ref{kernel}). We first prove a weak type estimate~.

\begin{prop}\label{BWL} Let $M$ be a complete Riemannian manifold satisfying $(D)$ and admitting a Poincar\'e inequality $(P_1)$. Let $T$ be a linear operator which is bounded from $\dot{W}^{1,2}$ (resp. $W^{1,2}$) to $L^2$ and is associated to a kernel satisfying (\ref{kernel}). \\
Then $T$ is bounded from $\dot{W}^{1,1}$ (resp. $W^{1,1}$) to $L^{1,\infty}$.
\end{prop}
\dem We give the proof in the homogeneous case, it is the same in the non-homogeneous case. Let $f\in \dot{W}^{1,1}$. We want to show that
$$
\mu (\left\{ x\in M; |Tf(x)|>\alpha\right\})\lesssim \frac{1}{\alpha}\|\nabla f\|_{L^1} .
$$
 Take the Calder\'on-Zygmund decomposition -- homogeneous version of Proposition \ref{CZ} -- of $f$ for $\alpha>0$. We have
 $$Tf=Tg+T(\sum_ib_i) $$
 and $\left\lbrace |Tf|>\alpha\right\rbrace\subset \left\lbrace |Tg|>\frac{\alpha}{2}\right\rbrace\bigcup \left\lbrace |T(\sum_ib_i)|>\frac{\alpha}{2}\right\rbrace$.
 Since $T$ is bounded from $\dot{W}^{1,2}$ to $L^2$ then
 $$\mu\left(\left\lbrace |Tg|>\frac{\alpha}{2}\right\rbrace\right)\leq \frac{4}{\alpha^2}\int_{M}|Tg|^2d\mu\lesssim \frac{1}{\alpha^2}\|T\|_{\dot{W}^{1,2}\rightarrow L^{2}}\alpha \|\nabla f\|_{L^1}.$$
For $|T(\sum_ib_i)|=|\sum_{i}Tb_i|\leq \sum_i|Tb_i|$ we have
 \begin{align*}
 \mu\left(\left\lbrace|T(\sum_ib_i)|>\frac{\alpha}{2}\right\rbrace\right)&\leq \mu\left(\left\lbrace\sum_i|Tb_i|>\frac{\alpha}{2}\right\rbrace\right)
 \\
 &\leq\mu(\bigcup_i 4Q_i)+\mu\left(\{(M\backslash \cup_i 4Q_i);\sum_i|Tb_i|>\frac{\alpha}{2}\}\right).
 \end{align*}
  From $(D)$ and the homogeneous analog of (\ref{sbaze}) of Proposition \ref{CZ}, we have $\mu(\bigcup_i 4Q_i)\leq \frac{C}{\alpha}\|\nabla f\|_1$.
 It remains to estimate $\mu(A)=\mu(\left\lbrace (M\backslash \cup_i4Q_i);\,\sum_i|Tb_i|>\frac{\alpha}{2}\right\rbrace).$
We have \begin{equation*} A\subset \left\lbrace \sum_i \ind_{M\backslash 4Q_i}|Tb_i|>\frac{\alpha}{2}\right\rbrace.
 \end{equation*}
 Then
 $$\mu(A)\leq \frac{2}{\alpha}\int_{M}\sum_i |Tb_i|\ind_{M\backslash4Q_i}d\mu=\frac{2}{\alpha}\sum_i\int_{M\backslash4Q_{i}}|Tb_i|d\mu.
 $$
Let $y_i\in Q_i$ such that $K(x,y_i)$ exists. Noting that $\int b_id\mu=0$, it comes that
 \begin{align*}
 \int_{M\backslash4Q_{i}}|Tb_i|(x)d\mu(x)&=\int_{M\backslash4Q_{i}}\left|\int_{Q_i}K(x,y)b_i(y)d\mu(y)\right|d\mu(x)
 \\
 &=\int_{M\backslash4Q_{i}}\left|\int_{Q_{i}}(K(x,y)-K(x,y_i))b_i(y)d\mu(y)\right|d\mu(x)
 \\
 &\leq \int_{Q_i}\left(\int_{M\backslash 4Q_i}|K(x,y)-K(x,y_i)|d\mu(x)\right)|b_i(y)|d\mu(y)
 \\
 &\lesssim \frac{1}{r_i}\int_{Q_i}|b_i(y)|d\mu(y)\sup_{y,\,y_i\in Q_i}r_i\int_{M\backslash4Q_i}|K(x,y)-K(x,y_i)|d\mu(x)
 \\
 &\lesssim \alpha \mu(Q_i).
 \end{align*}
Summing over $i$ and using the homogeneous analogous property of (\ref{sbaze}), the  proof is therefore complete. \findem

\mb To obtain this weak type estimate, we have to assume a strong Poincar\'e inequality $(P_1)$. The result of Theorem
\ref{HSL}  is also interesting: we are able to obtain a strong type estimate using Hardy-Sobolev spaces (instead of
the Sobolev space $\dot{W}^{1,1}$), and requiring a weaker Poincar\'e inequality in the non-homogeneous case.

\mb \textbf{Proof of Theorem \ref{HSL}:}
 We begin showing that in both case item 1. (resp. 2.), there exists a constant $C$, such that for all $2$-homogeneous atom $a$ (resp. non-homogenous atom),
\begin{equation}\label{batom}
\|Ta\|_{L^1}\leq C.
\end{equation}
We give the proof in the homogeneous case, it works the same in the non-homogeneous case.
Indeed, noting $Q=Q(x_0,r)$ the ball associated to the $(1,2)$ homogeneous atom $a$, we have
\begin{equation*}
\int_{4Q}|Ta|d\mu\leq C\|T\|_{\dot{W}^{1,2}\rightarrow L^2}\|a\|_{\dot{W}^{1,2}}{\mu(Q)^\frac{1}{2}}\leq C\|T\|_{\dot{W}^{1,2}\rightarrow L^2}.
\end{equation*}
On $M\backslash4Q$, we use the integral representation. The fact that $\int ad\mu=0$ yields
\begin{align*}
\int_{M\backslash4Q}|Ta|d\mu &\leq \int_{M\backslash4Q}\left|\int_{Q}K(x,y)a(y)d\mu(y)\right|d\mu(x)
 \\
 &=\int_{M\backslash4Q}\left|\int_{Q}(K(x,y)-K(x,x_0))a(y)d\mu(y)\right|d\mu(x)
 \\
 &\leq \int_{M\backslash4Q}\int_{Q}|K(x,y)-K(x,x_0)|\,|a(y)|d\mu(y)d\mu(x)
 \\
 &=\int_{Q}|a(y)|\left(\int_{M\backslash4Q}|K(x,y)-K(x,x_0)|d\mu(x)\right)d\mu(y)
 \\
 &= \int_{Q}|a-a_{Q}|\left(\int_{M\backslash4Q}|K(x,y)-K(x,x_0)|d\mu(x)\right)d\mu(y)
 \\
 &\leq Cr\mu(Q)(\aver{Q}|\nabla a|^{2}d\mu)^{\frac{1}{2}}\frac{C}{r}
 \\
 &\leq C.
 \end{align*}
 We used Poincar\'e inequality $(P_2)$, (\ref{kernel}) and the definition of a $(1,2)$ atom.

 \mb Now we conclude the proof of item 1. \\
%
 Thanks to Proposition \ref{BWL}, $T$ is bounded from $\dot{W}^{1,1}$ to $L^{1,\infty}$. Take  $f\in \dot{HS}_{(2),ato}^1:\,f=\sum_{i=1}^{\infty}\lambda_ib_i$ with for each $i$, $b_i$ is a $(1,2)$ homogeneous atom and with $\sum_{i=1}^{\infty}|\lambda_i|\sim\|f\|_{\dot{HS}_{(2),ato}^1}$. Since $\dot{HS}_{(2),ato}^1 \hookrightarrow \dot{W}^{1,1}$, we know that  $f_N=\sum_{i=1}^{N}\lambda_ib_i\in\dot{HS}_{F,(2),ato}^1$ converges to $f$ in $\dot{W}^{1,1}$. Thus by Proposition \ref{BWL},  $Tf_{N}$ converges to $Tf$ in $L^{1,\infty}$.

 \mb On the other hand, $Tf_{N}$ converges to $\sum_{i=1}^{\infty}\lambda_i Tb_i$ in $\dot{W}^{1,1}$ and therefore $Tf= \sum_{i=1}^{\infty}\lambda_i Tb_i$ and $\|Tf\|_{1}\leq C\|f\|_{\dot{HS}_{(2),ato}^1}$.

 \mb It remains to complete the proof of item 2.
 For this, we invoke the following lemma which finishes the proof. It is a Sobolev version of a result in \cite{MSV}, that was generalized in \cite{B2}.
 \findem

\begin{lem} \label{lS} Assume that $(P_2)$ holds.
Let $T$ be a bounded linear operator  from $W^{1,2}$ to $L^2$ with a constant $C$ such that for all $(1,2)$ atom
$f\in HS_{F,(2),ato}^1$, we have
 $$ \|T(f)\|_{L^1} \leq C.$$
Then $T$ extends continuously from $HS_{(2),ato}^1$ into $L^1$.
\end{lem}
\begin{rem} The proof uses the embedding $HS_{(2),ato}^1 \hookrightarrow L^1$, which does not hold for the homogeneous
space $\dot{HS}_{(2),ato}^1$. Actually, we do not know if such a result is true or not for homogeneous Hardy-Sobolev
spaces, without using (as it is well-known) a weak-type inequality from $\dot{W}^{1,1}$ to $L^{1,\infty}$ which
requires the Poincar\'e inequality $(P_1)$ as we saw in item 1.
\end{rem}

\dem  As $HS_{F,(2),ato}^1$ is dense in $HS_{(2),ato}^1$, we know that there exists an operator $U$ bounded from $HS_{(2),ato}^1$ into $L^1$ such that for each atom $m$: $U(m)=T(m)$. We have to prove that
$$ \forall f\in W^{1,2} \cap HS_{(2),ato}^1, \qquad U(f)=T(f).$$
To prove this fact, we use duality. Let $Q$ be a ball and $\phi_Q$ be a smooth function supported in $Q$ verifying
$$ \int_Q \phi_Q d\mu =1, \qquad \|\phi_Q\|_\infty\lesssim \frac{1}{\mu(Q)}, \qquad  \|\nabla \phi_Q\|_{\infty}\lesssim \frac{1}{r_Q\mu(Q)}.$$
Then for all smooth function $k$ supported in $Q$, with $\|k\|_{W^{1,2}}\leq \mu(Q)^{-1/2}$, the function $h:=k-(\int_Q
k)\phi_Q$ is a $(1,2)$-atom associated to the ball $Q$ (due to Poincar\'e inequality and Proposition \ref{prop:moy}). Let
$g\in L^\infty \cap L^2$. We have
$$ \langle T(h),g\rangle = \langle U(h), g\rangle.$$
We deduce that
$$ \langle h,T^*g\rangle = \langle h,U^* g\rangle.$$
Hence
$$ \left\langle k,\left[T^*g-U^*g\right] - \left(\int\phi_Q \left[T^*g-U^*g\right] d\mu\right){\bf 1}_Q \right \rangle =0.$$
We set $\lambda$ for the function $\lambda:=\left[T^*g-U^*g\right]$. We have 
$$ \left\|\lambda - \left(\int\phi_Q \lambda d\mu\right){\bf 1}_Q \right\|_{W^{-1,2}(B)}=0.$$
Thus $\lambda$ (as distribution) is constant on the ball $Q$. This fact is proved for every ball $Q$. We conclude that $\lambda$ (which is independent with respect to the ball) is constant over the whole manifold $M$. \\
The non-homogeneous Hardy-Sobolev space $HS_{(2),ato}^1$ is embedded into $L^1$. Then by $L^1-L^\infty$ duality, for
all functions $h\in HS_{(2),ato}^1$ we have
$$ \langle h, \lambda \rangle=0.$$
In particular for $f\in W^{1,2} \cap HS_{(2),ato}^1$, we get
\begin{align*} 
\langle f, \lambda \rangle & = 0 = _{\dot{W}^{1,2}}\langle f,T^*g\rangle_{\dot{W}^{-1,2}} - _{HS_{(2),ato}^1}\langle f,U^*g\rangle_{(HS_{(2),ato}^1)^*} \\
& = _{L^{2}}\langle T(f),g\rangle_{L^{2}} - _{L^1}\langle U(f),g\rangle_{L^\infty}.
\end{align*}
This is true for all functions $g\in L^\infty \cap L^{2}$. We deduce that $T(f)=U(f)$ in $\left(L^\infty \cap
L^{2}\right)^*$ and therefore $T(f)(x)=U(f)(x)$ for almost every $x\in M$. \findem

\gb \textbf{Proof of Corollary \ref{cor1}:} The proof follows from  the interpolation results in Theorem
\ref{thm:hardy} and Theorem \ref{thm:hardyhomo} and the self-improvement of Poincar\'e inequality of Theorem \ref{kz}.
\findem

\mb The result  in item 2. of Corollary \ref{cor1} can  also be  recovered  by suitably choosing the  operators $B_Q$
of the abstract Hardy-Sobolev spaces (defined in Subsection \ref{subsec:method1}).

\begin{df} For each ball $Q$ of $M$, we define our operator $B_Q$ as~:
$$ B_Q(f):=f- \left(\int_Q f d\mu \right) \phi_Q, $$
where $\phi_Q$ is a smooth function supported in $Q$ such that
$$ \int_Q \phi_Q d\mu =1, \qquad \|\phi_Q\|_\infty\lesssim \frac{1}{\mu(Q)}, \qquad  \|\nabla \phi_Q\|_{\infty}\lesssim r_Q^{-1}\mu(Q)^{-1}.$$
\end{df}

\mb With $\B=2$, we define our Hardy-Sobolev space $HW^1_{(2),ato}$. \\
We check the desired assumptions. Thanks to the Proposition \ref{prop:moy}, it is clear that under Poincar\'e inequality $(P_{2})$ the operators $A_Q$ are uniformly bounded on $W^{1,2}$. \\
Then by similar arguments as that in the proof of Theorem \ref{HSL}, under $(P_2)$ the above operator $T$ admits a continuous extension from $HW^1_{(2),ato}$ to $L^1$. Moreover, for $q\in(1,2)$ the inequality $(P_{q})$ implies that the maximal operator $M_{{\mathbb B},q'}$ is bounded by $\mathcal{M}_{2}$ (using Proposition \ref{prop:utile}). Using Theorem \ref{thm:generalbis}, we recover item 2. of Corollary \ref{cor1}.

\subsection{Application: $(RR_p)$.}
Let $M$ be a complete Riemannian manifold satisfying $(D)$.  Consider the linear operator $\Delta^{\frac{1}{2}}$ with
the following resolution
$$
\Delta^{\frac{1}{2}}f= c \int_{0}^{\infty}\Delta e^{-t\Delta}f\frac{dt}{\sqrt{t}} ,\quad  f\in C^{\infty}_{0}
$$
where $c=\pi^{-\frac{1}{2}}$. Here $\Delta^{\frac{1}{2}}f$ can be defined for $f\in \Lip$ as a measurable function (see
\cite{AC}). Since $\Delta^{\frac{1}{2}}1=0$, $\Delta^{\frac{1}{2}}$ can be defined on $\Lip\cap\dot{W}^{1,q}$ by taking
quotient which we keep calling $\Delta^{\frac{1}{2}}$. Applying Theorem \ref{HSL}, we obtain the following theorem for
$\Delta^{\frac{1}{2}}$.
\begin{thm}1- Let $M$ be a complete Riemannian manifold satisfying $(D)$ and $(P_1)$. Then  $\Delta^{\frac{1}{2}}$ is bounded from $\dot{HS}_{(r),ato}^{1}$ to $L^1$ for any  $r>1$.
\\
Consequently, $\Delta^\frac{1}{2}$ is bounded from $\dot{W}^{1,p}$ to $L^p$ for any $p\in(1,2]$.
\\
2-Let $M$ be a complete Riemannian manifold satisfying $(D)$ and $(P_q)$ for some $q\in[1,2)$. Then  $(I+\Delta)^{\frac{1}{2}}$ is bounded from $HS_{(r),ato}^{1}$ to $L^1$ for any $r\geq q$ if $q\neq 1$ (resp. $r>1$ if $q=1$).
\\
Consequently, $(I+\Delta)^\frac{1}{2}$ is bounded from $W^{1,p}$ to $L^p$ for any $p\in[q,2]$.
\end{thm}

\begin{rem} We refer the reader to \cite{ACDH,AC} for the study of inequality $(RR_p)$ for $p\in(1,2]$ (which corresponds to the boundedness of $\Delta^\frac{1}{2}$ from $\dot{W}^{1,p}$ to $L^p$) under Poincar\'e
inequality. The new point here is the limit case $(RR_1)$.
\end{rem}

\dem We prove item 1. of this theorem. We proceed analogously for the proof of item 2. Let us check that $\Delta^{\frac{1}{2}}$ satisfies the hypotheses of Theorem \ref{HSL}. First $\Delta^{\frac{1}{2}}$ is bounded from $\dot{W}_{2}^{1}$ to $L^2$. The kernel of $\Delta^{\frac{1}{2}}$ is $\int_{0}^{\infty}\partial_{t}p_{t}(x,y)\frac{dt}{\sqrt{t}}$. Under our hypotheses, the partial derivative of the heat kernel $\partial_{t}p_t$ verifies
\begin{equation}\label{kD}
|\partial_{t}p_{t}(x,y)|\leq \frac{C}{t\mu(B(y,\sqrt{t}))}e^{-\alpha\frac{d^{2}(x,y)}{t}}
\end{equation}
for every $x,\, y\in M$ and $t>0$ (see \cite{Davies2}, Theorem 4 and \cite{gri}, Corollary 3.3). Let $Q$ a ball of radius $r>0$ and $y,\,z\in Q$.  We therefore have
\begin{align*}
\int_{M\backslash4Q}\left|\int_{0}^{\infty}\partial_{t}(p_t(x,y)-p_t(x,z))\frac{dt}{\sqrt{t}}dt\right|d\mu(x)
& \\
& \hspace{-2cm} \leq\int_{M\backslash4Q}\int_{0}^{\infty}|\partial_{t}(p_t(x,y)-p_t(x,z))|\frac{dt}{\sqrt{t}}dtd\mu(x)
\\
&\hspace{-2cm} \leq C\int_{M\backslash4Q}\int_{0}^{\infty}\frac{1}{t\mu(Q(y,\sqrt{t}))}e^{-\alpha \frac{d^2(x,y)}{t}}\frac{dt}{\sqrt{t}}dtd\mu(x)
\\
&\hspace{-2.5cm} +C\int_{M\backslash4Q}\int_{0}^{\infty}\frac{1}{t\mu(Q(z,\sqrt{t}))}e^{-\alpha \frac{d^2(x,z)}{t}}\frac{dt}{\sqrt{t}}dtd\mu(x).
\end{align*}
Let us estimate  $I=\int_{M\backslash4Q}\left(\int_{0}^{\infty}\frac{1}{t\mu(Q(y,\sqrt{t}))}e^{-\alpha \frac{d^2(x,y)}{t}}\frac{dt}{\sqrt{t}}\right)d\mu(x)$. Since $y\in Q$ and $x\in M\backslash4Q$ then $d(x,y)\geq 3r$. It follows that
\begin{align*}
I&\leq \int_{0}^{\infty}\frac{C}{t\mu(Q(y,\sqrt{t}))}\left(\int_{\{x;\,d(x,y)>\sqrt{9r^2}\}}e^{-\alpha\frac{d^2(x,y)}{t}}d\mu(x) \right)\frac{dt}{\sqrt{t}}
\\
&\leq \int_{0}^{\infty}\frac{C}{t\mu(Q(y,\sqrt{t}))}C_{\alpha}\mu(Q(y,\sqrt{t}))e^{-\alpha\frac{9r^2}{t}}\frac{dt}{\sqrt{t}}
\\
&\leq \int_{0}^{\infty}\frac{e^{-\alpha\frac{9r^2}{t}}}{t\sqrt{t}}dt
\\
&\leq \frac{C}{r}\int_{0}^{\infty}\frac{e^{-\alpha \frac{9}{t}}}{t\sqrt{t}}dt
\\
&\leq \frac{C}{r}.
\end{align*}
In the second estimate, we used that $\int_{d(x,y)>\sqrt{t}}e^{-\gamma \frac{d^{2}(x,y)}{s}}d\mu(x)\leq
C_{\gamma}\mu(Q(y,\sqrt{s}))e^{-\gamma\frac{t}{s}}$ (\cite{CD1}, Lemma 2.1 ). Similarly, we prove that
$\int_{M\backslash4Q}\left(\int_{0}^{\infty}\frac{1}{t\mu(Q(z,\sqrt{t}))}e^{-\alpha
\frac{d^2(x,z)}{t}}\frac{dt}{\sqrt{t}}\right)d\mu(x)\leq \frac{C}{r}$. Taking the supremum over all $y,z\in Q$, all
balls $Q$ and applying Theorem \ref{HSL}, we obtain that $T$ is bounded from $\dot{HS}_{r,ato}^{1}$ to $L^1$ for $r>1$.  Finally the boundedness of $T$ from $\dot{W}^{1,p}$ to $L^p$ for
$1<p<2$ follows from Corollary \ref{cor1}.
 \findem

\end{document}